\documentclass[a4paper,12pt]{article} 
\usepackage{amssymb}
\usepackage{amscd}
\usepackage{amsmath} 
\usepackage{graphicx}
\usepackage{latexsym} 
\usepackage{enumerate}

\usepackage{theorem}
\theoremstyle{plain}
\theorembodyfont{\itshape}
\newtheorem{theorem}{Theorem}

\newtheorem{lemma}[theorem]{Lemma}

\theorembodyfont{\rmfamily}

\newtheorem{remark}[theorem]{Remark}
\newtheorem{problem}[theorem]{Problem}
\newtheorem{example}[theorem]{Example}

\def\rP{\mathrm{P}}
\def\PVI{\mathrm{P}_{\mathrm{VI}}}
\def\rQ{\mathrm{Q}}
\def\rR{\mathrm{R}}

\def\RH{\mathrm{RH}} 
\def\rh{\mathrm{rh}} 

\def\E{\mathcal{E}}

\def\I{\mathcal{I}}
\def\K{\mathcal{K}}
\def\M{\mathcal{M}}
\def\O{\mathcal{O}}

\def\Sol{\mathcal{S}}
\def\C{\mathbb{C}}
\def\Q{\mathbb{Q}}

\def\P{\mathbb{P}}
\def\bR{\mathbb{R}}

\def\Wall{\mathbf{Wall}} 
\def\0{\mathbf{0}}
\def\a{\alpha}

\def\ab{\mbox{\boldmath$a$}}
\def\bb{\mbox{\boldmath$b$}}
\def\cb{\mbox{\boldmath$c$}}

\def\ga{\gamma}
\def\k{\kappa}
\def\kb{\mbox{\boldmath$\kappa$}}

\def\si{\sigma}
\def\th{\theta}
\def\vD{\varDelta}

\def\Th{\Theta}
\def\ve{\varepsilon}
\def\carl{\circlearrowleft}
\def\car{\curvearrowright}

\def\llra{\longleftrightarrow}

\def\ol{\overline}
\def\ora{\overrightarrow}

\def\dfrac#1#2{{\displaystyle\frac{#1}{#2}}}

\def\bm{\boldmath}

\def\la{\langle}
\def\ra{\rangle}

\def\rot{\rotatebox}

\def\wt{\widetilde} 
\title{\bf On Algebraic Solutions to 
Painlev\'e VI}
\author{Katsunori Iwasaki \\ \\ 
Faculty of Mathematics, Kyushu University \\ 
6-10-1 Hakozaki, Higashi-ku, Fukuoka 812-8581 Japan} 
\date{}
\begin{document}
\maketitle
\begin{abstract} 
We announce some results which might bring a new insight into 
the classification of algebraic solutions to the sixth 
Painlev\'{e} equation. 
They consist of the rationality of parameters, 
trigonometric Diophantine conditions, and what the author 
calls the Tetrahedral Theorem regarding the absence of algebraic 
solutions in certain situations. 
The method is based on fruitful interactions between the moduli 
theoretical formulation of Painlev\'e VI and dynamics on 
character varieties via the Riemann-Hilbert correspondence. 
\end{abstract} 
\section{Introduction} \label{sec:intro}
All algebraic solutions to the Gauss hypergeometric equation 
were classified by H.A. Schwarz \cite{Schwarz} in 1873. 
After him this classifiaction has been known as Schwarz's list. 
On the other hand the sixth Painlev\'e equation is known 
as a nonlinear generalization of the Gauss equation. 
So we are naturally led to the problem of classifying 
all glgebraic solutions to Painlev\'e VI. 
This problem is still open (as of this writing) and 
there is a vast literature on this theme including 
\cite{AK,BeG,Boalch1,Boalch2,Boalch4,Doran,Dubrovin,DM,Hitchin1,
Hitchin2,Kitaev1,Kitaev2,Mazzocco1,YL}. 
The attempt at solving this problem could be entitled {\sl Towards a 
nonlinear Schwarz's list} as P. Boalch employs these words as 
the title of his survey \cite{Boalch3}, in which the current 
states of the subject are nicely presented. 
The aim of this article is to announce some new results which might 
bring a new insight into this subject. 
\par
The above-mentioned problem for Painlev\'e VI is closely related to 
a problem from topology, that is, to classifying all finite orbits of 
the mapping class group action on certain character varieties, where 
the Painlev\'e-equation side and the character-variety side are connected 
by the so-called Riemann-Hilbert correspondence. 
Our philosophy is that working on both sides together, going back 
and forth between them, should be more fruitful than working on 
only one side of them. 
The mixture of methods from both sides should go much farther 
than either side could go by itself. 
The main results of this article are the rationality of parameters 
(\S\ref{sec:rationality}), trigonometric Diophantine conditions 
(\S\ref{sec:tdc}), and what the author calls the 
Tetrahedral Theorem (\S\ref{sec:tetra}) which is concerned with 
the absence of algebraic solutions in certain situations. 
\par
The contents of this article are based on the following talks by 
the author: (1) a series of talks at IRMAR, l'Universit\'e de 
Rennes, March, 2008. 
The author thanks S. Cantat and F. Loray for stimulating discussions; 
(2) a talk at the Conference on Exact WKB Analysis and Microlocal 
Analysis in RIMS, Kyoto, May, 2008. 
This article is a contribution to its Proceedings; 
(3) a talk at the International Conference ``From Painlev\'e 
to Okamoto" in The University of Tokyo, June, 2008. 
A full account of this announcement will be given in \cite{Iwasaki2}. 
\par\medskip
{\bf Note}. After an earlier version \cite{IwasakiVer1} of this 
article had been posted in the e-Print arXiv, a preprint 
\cite{LT} giving a complete classification of algebraic 
solutions was posted in the arXiv by O.~Lisovyy and Y.~Tykhyy. 
Their approach is quite straightforward. 
First, they use trigonometric Diophantine conditions to show that 
all monodromy data that can lead to finite orbits necessarily belong 
to an explicitly defined finite set (with two exceptions), and then 
all possibilities are checked by computers. 
Hearing of their work, the author puts this note here instead of 
revising the Introduction to a large extent. 
\section{Dynamics on Character Varieties} \label{sec:char}
Let $X$ be a real orientable closed surface with a finite number 
of punctures. 
By definition a relative $SL_2(\C)$-character variety of $X$ is 
the moduli space of Jordan equivalence classes of representations 
into $SL_2(\C)$ of the fundamental group $\pi_1(X)$ with prescribed 
local representations around the punctures. 
Hereafter a relative $SL_2(\C)$-character variety is simply 
referred to as a character variety. 
It is acted on by the 
mapping class group of $X$ in a natural manner. 
\par 
In this article we are interested in the basic case where 
$X$ is the quadruply-punctured sphere. 
In this case the character varieties are realized as the 
four-parameter family of complex affine cubic surfaces 
$\Sol(\th) = \{\, x = (x_1,x_2,x_3) \in \C_x^3 
\,:\, f(x,\th) = 0 \, \}$ parametrized by 
$\th = (\th_1,\th_2,\th_3,\th_4) \in \Th := 
\C_{\th}^4$, where the polynomial $f(x,\th)$ is defined by 
\[
f(x,\th) := x_1x_2x_3+x_1^2+x_2^2+x_3^2
-\th_1x_1-\th_2x_2-\th_3x_3+\th_4. 
\]
The surface $\Sol(\th)$ is a $(2,2,2)$-surface, that is, 
the defininig function $f(x,\th)$ is a quadratic polynomial 
in each variable $x_i$ $(i = 1,2,3)$. 
Thus the line through a point $x \in \Sol(\th)$ parallel to 
the $x_i$-axis passes through a unique second point 
$x' = \si_i(x) \in \Sol(\th)$. 
This induces an involutive automorphism 
$\si_i : \Sol(\th) \to \Sol(\th)$ for each $i = 1,2,3$ 
(see Figure \ref{fig:involution}). 
\begin{figure}[t]
\begin{center}
\unitlength 0.1in
\begin{picture}(33.24,17.19)(4.90,-19.59)
%
\special{pn 13}%
\special{pa 490 1014}%
\special{pa 1111 1014}%
\special{fp}%
%
\special{pn 13}%
\special{pa 2641 1023}%
\special{pa 3397 1023}%
\special{fp}%
%
\special{pn 20}%
\special{ar 1858 1095 1053 855  3.3124258 6.1241532}%
%
\special{pn 20}%
\special{ar 1849 1104 1062 855  6.2694533 6.2831853}%
\special{ar 1849 1104 1062 855  0.0000000 3.1735968}%
\put(17.1400,-16.6200){\makebox(0,0)[lb]{$\Sol(\th)$}}%
%
\special{pn 13}%
\special{pa 1201 1014}%
\special{pa 2524 1023}%
\special{dt 0.045}%
\special{pa 2524 1023}%
\special{pa 2523 1023}%
\special{dt 0.045}%
\put(10.6000,-9.0000){\makebox(0,0)[lb]{$x$}}%
\put(25.8000,-9.1000){\makebox(0,0)[lb]{$x'$}}%
%
\special{pn 20}%
\special{sh 0.600}%
\special{ar 1102 1005 44 44  0.0000000 6.2831853}%
%
\special{pn 20}%
\special{sh 0.600}%
\special{ar 2632 1014 44 44  0.0000000 6.2831853}%
\put(18.1000,-6.1000){\makebox(0,0)[lb]{$\si_i$}}%
%
\special{pn 13}%
\special{pa 3040 1295}%
\special{pa 3814 1295}%
\special{fp}%
\special{sh 1}%
\special{pa 3814 1295}%
\special{pa 3747 1275}%
\special{pa 3761 1295}%
\special{pa 3747 1315}%
\special{pa 3814 1295}%
\special{fp}%
\put(31.8400,-15.2900){\makebox(0,0)[lb]{$x_i$-axis}}%
%
\special{pn 13}%
\special{ar 1880 1850 1161 1161  4.1675251 5.2764727}%
%
\special{pn 13}%
\special{pa 2380 800}%
\special{pa 2510 870}%
\special{fp}%
\special{sh 1}%
\special{pa 2510 870}%
\special{pa 2461 821}%
\special{pa 2463 845}%
\special{pa 2442 856}%
\special{pa 2510 870}%
\special{fp}%
\end{picture}%
\end{center}
\caption{Involutions on the $(2,2,2)$-surface $\Sol(\th)$}
\label{fig:involution}
\end{figure}
Let $G$ be the group generated by these three involutions. 
Then it turns out that the generators have no other relations 
than the trivial ones $\si_1^2 = \si_2^2 = \si_3^2 =1$. 
Namely we have
\[
G := \la \si_1, \si_2, \si_3 \ra = 
\la \si_1, \si_2, \si_3 \,|\, 
\si_1^2 = \si_2^2 = \si_3^2 =1 \ra \car \Sol(\th). 
\]
Each element $\si \in G$ can be written in a unique 
way as a word 
$\si = \si_{i_1}\si_{i_2}\cdots\si_{i_n}$ with alphabet 
$\{\si_1,\si_2,\si_3\}$ such that the consecutive indices 
$i_{\nu}$ and $i_{\nu+1}$ are all distinct. 
Let $G(2)$ denote the subgroup of all even words in $G$. 
It is an index-two normal subgroup of $G$. 
In the present case the mapping class group action is realized 
as the group action $G(2) \car \Sol(\th)$.  
\begin{problem} \label{prob1}
Classify all finite orbits of the action $G(2) \car \Sol(\th)$. 
\end{problem}
\par
Let $V = \{\, \th \in \Th \,:\, \vD(\th) = 0 \,\}$ be the 
discriminant locus of the family of cubic surfaces 
$\Sol(\th)$ parametrized by $\th \in \Th$, where $\vD(\th)$ 
is the discriminant of $f(x,\th)$ as a polynomial of $x$. 
For any $\th \in V$ the surface $\Sol(\th)$ has at most 
four simple singulatities. 
Let 
\begin{equation} \label{eqn:varphi}
\varphi : \wt{\Sol}(\th) \to \Sol(\th) 
\end{equation} 
be an algebraic minimal desingularization. 
Then the action $G \car \Sol(\th)$ lifts to the smooth surface 
$\wt{\Sol}(\th)$ in a unique way and Problem \ref{prob1} is 
refined into the following problem. 
\begin{problem} \label{prob2}
Classify all finite orbits of the lifted action 
$G(2) \car \wt{\Sol}(\th)$. 
\end{problem}
\par 
It is easy to see that the singular points of $\Sol(\th)$ 
are exactly the fixed points of the action $G(2) \car \Sol(\th)$ 
so that the exceptional set $\E(\th) \subset \wt{\Sol}(\th)$ 
is invariant under the lifted action $G(2) \car \wt{\Sol}(\th)$. 
Problem \ref{prob2} is finer than Problem \ref{prob1} to the 
extent that Problem \ref{prob2} demands to classify 
finite orbits on the exceptional set $\E(\th)$. 
But this extra task is not so heavy as will be explained in 
\S\ref{sec:RH}. 
So one can safely say that the two problems are 
approximately the same. 
\section{The Sixth Painlev\'e Equation} \label{sec:PVI}
The sixth Painlev\'e equation $\PVI(\k)$ is a non-autonomous 
Hamiltonian system with a complex time variable 
$z \in Z := \P^1-\{0,1,\infty\}$ and unknown functions 
$q = q(z)$ and $p = p(z)$, 
\[
\dfrac{dq}{dz} = \dfrac{\partial H(\k)}{\partial p}, \qquad 
\dfrac{dp}{dz} = -\dfrac{\partial H(\k)}{\partial q}, 
\]
depending on complex parameters $\k$ in the $4$-dimensional 
affine space 
\[
\K := \{\,\k = (\k_0,\k_1,\k_2,\k_3,\k_4) \in \C_{\k}^5 \,:\, 
2 \k_0 + \k_1 + \k_2 + \k_3 + \k_4 = 1 \,\}, 
\]
where the Hamiltonian $H(\k) = H(q,p,z;\k)$ is given by 
\[
z(z-1) H(\k) = (q_0q_zq_1)p^2-
\{\k_1 q_1q_z +(\k_2-1) q_0q_1 + \k_3 q_0 q_z \}p + 
\k_0(\k_0+\k_4) q_z
\]
with $q_{\nu} := q - \nu$ for $\nu \in \{0,z,1\}$. 
It is known that $\PVI(\k)$ has the Painlev\'e property in 
$Z$, that is, any meromorphic solution germ to $\PVI(\k)$ at 
a base point $z \in Z$ admits a global analytic continuation 
along any path in $Z$ emanating from $z$ as a meromorphic 
function. 
In fact, this property is a natural consequence of our 
solution to the Riemann-Hilbert problem based on a suitable 
moduli theoretical formulation of the sixth Painlev\'e 
equation (see \cite{IIS1,IIS2}). 
\par
For the Painlev\'{e} equation we are interested in the 
following problem. 
\begin{problem} \label{prob3}
Classify all algebraic solutions to $\PVI(\k)$. 
\end{problem}
For the current state of the problem we refer to the 
nice survey article \cite{Boalch3}. 
We also consider a closely related problem (which turns out to 
be an equivalent problem). 
Fix a base point $z \in Z$ and let $\M_z(\k)$ be the set of 
all meromorphic solution germs to $\PVI(\k)$ at the point $z$. 
Thanks to the Painlev\'e property, any germ $Q \in \M_z(\k)$ 
can be continued analytically along any loop $\ga \in \pi_1(Z,z)$ 
into a second germ $\ga_*Q \in \M_z(\k)$. 
This defines an automorphism $\ga_* : \M_z(\k) \carl$ and 
hence a group action $\pi_1(Z,z) \car \M_z(\k)$, called the 
nonlinear monodromy action. 
\begin{problem} \label{prob4} 
Classify all finite orbits of the nonlinear monodromy action 
$\pi_1(Z,z) \car \M_z(\k)$. 
\end{problem}
\par
Since any algebraic solution to $\PVI(\k)$ has only finitely 
many local branches at the base point $z$ and these local 
branches are permuted by the $\pi_1(Z,z)$-action, there is 
the natural inclusion: 
\begin{equation} \label{eqn:incl}
\{\, \mbox{germs at $z$ of algebraic solutions to $\PVI(\k)$}\,\} 
\hookrightarrow 
\{\, \mbox{finite $\pi_1(Z,z)$-orbits on $\M_z(\k)$}\,\} 
\end{equation}
One may be worried about the difference of the two sets. 
In fact there is no difference. 
\begin{theorem}[\cite{Iwasaki1}] \label{thm:alg=finite}
The inclusion $(\ref{eqn:incl})$ is surjective and hence 
Problems $\ref{prob3}$ and $\ref{prob4}$ are equivalent. 
\end{theorem}
There is a small gap in an argument of \cite{Iwasaki1}, 
which is to be filled in \cite{Iwasaki2}. 
\section{Riemann-Hilbert correspondence} \label{sec:RH}
To connect Problem \ref{prob2} with Problem \ref{prob3} (or 
equivalently with Problem \ref{prob4}), we review the 
Riemann-Hilbert correspondence \cite{IIS1,IIS2,Iwasaki1}. 
It exists in the parameter level and in the moduli level. 
\par
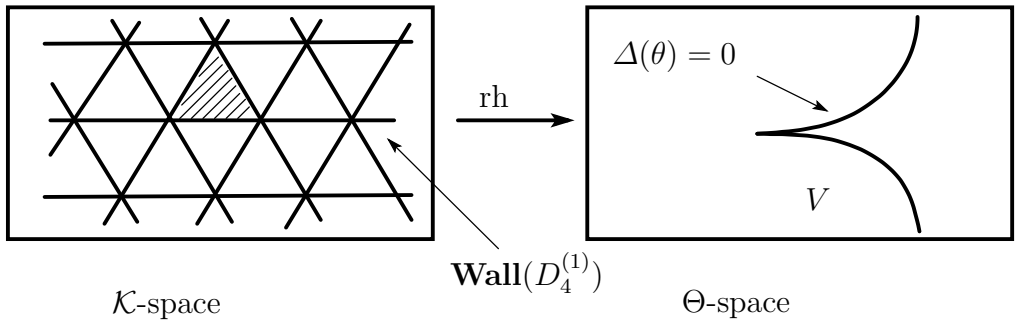
\begin{figure}[t]
\begin{center}
\unitlength 0.1in
\begin{picture}( 52.0000, 14.6000)(  2.0000,-20.7000)
%
\special{pn 20}%
\special{pa 200 610}%
\special{pa 2400 610}%
\special{pa 2400 1820}%
\special{pa 200 1820}%
\special{pa 200 610}%
\special{fp}%
%
\special{pn 20}%
\special{pa 3200 610}%
\special{pa 5400 610}%
\special{pa 5400 1820}%
\special{pa 3200 1820}%
\special{pa 3200 610}%
\special{fp}%
%
\special{pn 20}%
\special{pa 420 1200}%
\special{pa 2200 1200}%
\special{fp}%
%
\special{pn 20}%
\special{pa 420 1380}%
\special{pa 880 690}%
\special{fp}%
%
\special{pn 20}%
\special{pa 740 680}%
\special{pa 1360 1730}%
\special{fp}%
%
\special{pn 20}%
\special{pa 1340 680}%
\special{pa 700 1740}%
\special{fp}%
%
\special{pn 20}%
\special{pa 1200 670}%
\special{pa 1820 1720}%
\special{fp}%
%
\special{pn 20}%
\special{pa 1820 680}%
\special{pa 1180 1740}%
\special{fp}%
%
\special{pn 20}%
\special{pa 1670 670}%
\special{pa 2290 1720}%
\special{fp}%
%
\special{pn 20}%
\special{pa 2280 700}%
\special{pa 1640 1760}%
\special{fp}%
%
\special{pn 20}%
\special{pa 870 1730}%
\special{pa 430 1010}%
\special{fp}%
%
\special{pn 20}%
\special{pa 390 800}%
\special{pa 2290 790}%
\special{fp}%
%
\special{pn 20}%
\special{pa 390 1600}%
\special{pa 2290 1590}%
\special{fp}%
%
\special{pn 8}%
\special{pa 1380 1020}%
\special{pa 1210 1190}%
\special{fp}%
\special{pa 1360 980}%
\special{pa 1150 1190}%
\special{fp}%
\special{pa 1340 940}%
\special{pa 1090 1190}%
\special{fp}%
\special{pa 1310 910}%
\special{pa 1100 1120}%
\special{fp}%
\special{pa 1290 870}%
\special{pa 1190 970}%
\special{fp}%
\special{pa 1400 1060}%
\special{pa 1270 1190}%
\special{fp}%
\special{pa 1430 1090}%
\special{pa 1330 1190}%
\special{fp}%
\special{pa 1450 1130}%
\special{pa 1390 1190}%
\special{fp}%
\put(26.4000,-11.3000){\makebox(0,0)[lb]{$\rh$}}%
\put(33.3000,-9.5000){\makebox(0,0)[lb]{${\mit\Delta}(\theta) = 0$}}%
\put(7.4000,-22.4000){\makebox(0,0)[lb]{$\K$-space}}%
\put(36.9000,-22.3000){\makebox(0,0)[lb]{$\Theta$-space}}%
\put(24.9000,-21.0000){\makebox(0,0)[lb]{$\Wall(D_4^{(1)})$}}%
%
\special{pn 8}%
\special{pa 4080 970}%
\special{pa 4440 1150}%
\special{fp}%
\special{sh 1}%
\special{pa 4440 1150}%
\special{pa 4390 1102}%
\special{pa 4392 1126}%
\special{pa 4372 1138}%
\special{pa 4440 1150}%
\special{fp}%
%
\special{pn 20}%
\special{pa 4080 1270}%
\special{pa 4112 1268}%
\special{pa 4146 1266}%
\special{pa 4178 1264}%
\special{pa 4210 1260}%
\special{pa 4242 1258}%
\special{pa 4274 1254}%
\special{pa 4306 1250}%
\special{pa 4338 1244}%
\special{pa 4368 1238}%
\special{pa 4400 1232}%
\special{pa 4430 1224}%
\special{pa 4462 1214}%
\special{pa 4492 1204}%
\special{pa 4522 1192}%
\special{pa 4550 1180}%
\special{pa 4580 1166}%
\special{pa 4608 1150}%
\special{pa 4636 1132}%
\special{pa 4662 1114}%
\special{pa 4688 1094}%
\special{pa 4714 1074}%
\special{pa 4738 1052}%
\special{pa 4760 1028}%
\special{pa 4782 1004}%
\special{pa 4802 978}%
\special{pa 4820 952}%
\special{pa 4838 924}%
\special{pa 4854 896}%
\special{pa 4868 866}%
\special{pa 4880 836}%
\special{pa 4890 804}%
\special{pa 4898 772}%
\special{pa 4904 742}%
\special{pa 4908 710}%
\special{pa 4910 678}%
\special{pa 4910 660}%
\special{sp}%
%
\special{pn 20}%
\special{pa 4080 1270}%
\special{pa 4114 1270}%
\special{pa 4146 1272}%
\special{pa 4178 1272}%
\special{pa 4210 1272}%
\special{pa 4242 1274}%
\special{pa 4276 1276}%
\special{pa 4308 1278}%
\special{pa 4340 1280}%
\special{pa 4370 1284}%
\special{pa 4402 1288}%
\special{pa 4434 1294}%
\special{pa 4466 1300}%
\special{pa 4496 1308}%
\special{pa 4528 1318}%
\special{pa 4558 1328}%
\special{pa 4588 1340}%
\special{pa 4618 1352}%
\special{pa 4648 1368}%
\special{pa 4676 1384}%
\special{pa 4704 1402}%
\special{pa 4730 1420}%
\special{pa 4756 1442}%
\special{pa 4780 1464}%
\special{pa 4802 1488}%
\special{pa 4822 1512}%
\special{pa 4838 1538}%
\special{pa 4854 1566}%
\special{pa 4868 1594}%
\special{pa 4880 1624}%
\special{pa 4890 1654}%
\special{pa 4898 1686}%
\special{pa 4906 1716}%
\special{pa 4914 1748}%
\special{pa 4920 1780}%
\special{sp}%
%
\special{pn 8}%
\special{pa 2740 1870}%
\special{pa 2170 1300}%
\special{fp}%
\special{sh 1}%
\special{pa 2170 1300}%
\special{pa 2204 1362}%
\special{pa 2208 1338}%
\special{pa 2232 1334}%
\special{pa 2170 1300}%
\special{fp}%
%
\special{pn 20}%
\special{pa 2530 1200}%
\special{pa 3100 1200}%
\special{fp}%
\special{sh 1}%
\special{pa 3100 1200}%
\special{pa 3034 1180}%
\special{pa 3048 1200}%
\special{pa 3034 1220}%
\special{pa 3100 1200}%
\special{fp}%
\put(43.3000,-16.7000){\makebox(0,0)[lb]{$V$}}%
\end{picture}%
\end{center}
\caption{The Riemann-Hilbert correspondence in the parameter level}
\label{fig:rh}
\end{figure}

Firstly, the parameter space $\K$ is acted on by the affine Weyl 
group $W(D_4^{(1)})$ of type $D_4^{(1)}$ and the Riemann-Hilbert 
correspondence in the parameter level is a holomorphic map 
$\rh : \K \to \Th$ that is a branched $W(D_4^{(1)})$-covering 
ramifying along $\Wall(D_4^{(1)})$ and mapping it onto the 
discriminant locus $V \subset \Th$ of the family of cubic surfaces, 
where $\Wall(D_4^{(1)})$ is the union of all reflecting hyperplanes 
for the reflection group $W(D_4^{(1)})$ (see Figure \ref{fig:rh}). 
Secondly, developing a suitable moduli theory \cite{IIS1,IIS2} 
allows us to realize the set $\M_z(\k)$ as the moduli space of 
(certain) stable parabolic connections and thereby to provide it 
with the structure of a smooth quasi-projective rational surface. 
The Riemann-Hilbert correspondence (in the moduli level), 
\begin{equation} \label{eqn:RH}
\RH_{z,\k} : \M_z(\k) \to \Sol(\th), \quad Q \mapsto \rho, \qquad 
\mbox{with} \quad \th = \rh(\k), 
\end{equation}
is defined to be the holomorphic map sending each connection $Q$ to 
its monodromy representation $\rho$ up to Jordan equivalence. 
A fundamental fact for the map (\ref{eqn:RH}) is the following. 
\begin{theorem}[\cite{IIS1,IIS2}] \label{thm:RH}
The Riemann-Hilbert correspondence $(\ref{eqn:RH})$ is a proper 
surjective holomorphic map that yields an analytic minimal 
resolution of simple singularities. 
\end{theorem}
\par
By the minimality of the resolution, the Riemann-Hilbert 
correspondence (\ref{eqn:RH}) uniquely lifts to a biholomorphism 
$ \wt{\RH}_{z,\k} : \M_z(\k) \to \wt{\Sol}(\th)$ such that 
the following diagram is commutative: 
\[
\begin{CD}
\M_z(\k) @> \wt{\RH}_{z,\k} >> \wt{\Sol}(\th) \\
@|     @VV \varphi V  \\
\M_z(\k) @> \RH_{z,\k} >> \Sol(\th) 
\end{CD}
\]
The lifted Riemann-Hilbert correspondence $\wt{\RH}_{z,\k}$ gives 
a (strict) conjugacy between the nonlinear monodromy action 
$\pi_1(Z,z) \car \M_z(\k)$ and the mapping class group action 
$G(2) \car \wt{\Sol}(\th)$. 
In these circumstances the exceptional set $\E_z(\k) \subset 
\M_z(\k)$ of the resolution (\ref{eqn:RH}) just corresponds to 
the exceptional set $\E(\th) \subset \wt{\Sol}(\th)$ of the 
resolution (\ref{eqn:varphi}). 
We remark that $\E_z(\k)$ parametrizes the so-called Riccati 
solutions to $\PVI(\k)$, namely, those solutions which can be 
written in terms of the Riccati equations associated with 
Gauss hypergeometric equations (see \cite{IIS1}). 
\par
The lifted Riemann-Hilbert correspondence together with 
Theorem \ref{thm:alg=finite} yields the diagram: 
\[
\begin{array}{ccc}
\{\, \mbox{germs at $z$ of algebraic solutions to $\PVI(\k)$} \,\} 
& = & 
\{\, \mbox{finite $\pi_1(Z,z)$-orbits on $\M_z(\k)$}\, \} \\[2mm]
 & & \mbox{bijection} \,\, 
\rotatebox[origin=c]{90}{$\llra$} \,\, 
\wt{\RH}_{z,\k} \\[2mm]
 & & \{\, \mbox{finite $G(2)$-orbits on $\wt{\Sol}(\th)$}\, \}
\end{array}
\]
In summary, Problem \ref{prob1} is almost equivalent to 
Problem \ref{prob2}, while Problems \ref{prob2}, \ref{prob3} and 
\ref{prob4} are all equivalent. 
The difference of Problem \ref{prob2} from Problem \ref{prob1} 
amounts to classifying all Riccati algebraic solutions to 
$\PVI(\k)$, which in turn can be reduced to classifying all 
algebraic solutions to the Gauss hypergeometric equation, 
the classical problem settled by H.A.~Schwarz \cite{Schwarz}. 
\section{Rationality of Parameters} \label{sec:rationality}
An algebraic solution to $\PVI(\k)$ is said to be of degree 
$d$ if it has exactly $d$ local branches (germs) at a base 
point $z \in Z$. 
On the other hand a finite $G(2)$-orbit in $\wt{\Sol}(\th)$ is 
said to be of degree $d$ if it has exactly $d$ elements. 
Note that these two concepts of degree are consistent under 
the lifted Riemann-Hilbert correspondence 
$\wt{\RH}_{z,\k} : \M_z(\k) \to 
\wt{\Sol}(\th)$ with $\th = \rh(\k)$. 
\par
Naturally one may guess that those parameters $\k \in \K$ for 
which $\PVI(\k)$ admits at least one algebraic solutions of 
degree $d \ge d_0$ should have a very ``sparse" distribution, 
for some (perhaps reasonably large) integer $d_0$. 
Actually this guess is true in the following sense. 
\begin{theorem} \label{thm:rationality} 
If we take $d_0 = 7$, then we have the following 
rationality conditions. 
\begin{enumerate}
\item If $\PVI(\k)$ admits an algebraic solution of degree 
$d \ge 7$, then $\k_0$, $\k_1$, $\k_2$, 
$\k_3$ and $\k_4$ must be rational numbers.
\item If $\PVI(\k)$ admits an algebraic solution of degree 
$d \ge 1$ without univalent local branches at any fixed singular 
point $z = 0$, $1$, $\infty$, then $d\k_0$, $d\k_1$, $d\k_2$, 
$d\k_3$ and $d\k_4$ must be integers.
\end{enumerate}
\end{theorem}
Theorem \ref{thm:rationality} enables us to concentrate our 
attention on the rational and hence real part $\K_{\bR}$ of the 
complex affine space $\K$, as far as algebraic solutions of 
degree $d \ge 7$ are concerned. 
\begin{example} \label{ex:klein} 
To illustrate assertion (2) of Theorem \ref{thm:rationality}, 
we look at the ``Klein solution" constructed by 
Boalch \cite{Boalch1} based on the Klein complex reflection 
group of order $336$ in $SL_3(\C)$,  
\[
\left\{
\begin{array}{rcl}
z &=& \dfrac{(7s^2-7s+4)^2}{s^3(4s^2-7s+7)^2}, 
\\[6mm]
q &=& \dfrac{(s+1)(7s^2-7s+4)}{2s(s^2-s+1)(4s^2-7s+7)}, 
\\[6mm]
p &=& -\dfrac{2s(s+1)(s-2)(2s-1)(s^2-s+1)
(4s^2-7s+7)}{21(s-1)^2(4s^2-s+4)(7s^2-7s+4)}, 
\end{array}
\right.
\]
for which $d = 7$ and $\k=(1/7,1/7,1/7,1/7,2/7)$. 
This solution has ramification indices $3$, $2$, $2$ (a partition 
of $d = 7$) at each of the three fixed singular points 
$z = 0, 1, \infty$. 
Namely, it has one local branch of valency $3$ and 
two local branches of valency $2$ (and hence no 
univalent local branch) at each fixed singular point. 
Observe that $d \k_i$ $(i = 0,1,2,3,4)$ are 
integers. 
\end{example}
\par 
Two remarks are in order regarding Theorem \ref{thm:rationality}. 
\begin{remark} \label{rem:rationality}
For $i = 1$, $2$, item $(i)$ below is a remark about assertion 
$(i)$ of Theorem $\ref{thm:rationality}$. 
\begin{enumerate}
\item One may ask why condition $d \ge 7$ is imposed and 
how the assertion is derived. 
A brief answer to these questions will be given in 
\S\ref{sec:lines} (especially in Lemma \ref{lem:matrix} and 
the discussions thereafter). 
One may also ask what happens if $d \le 6$. 
It is known that there exist three exceptional classes of 
non-Riccati algebraic solutions to $\PVI(\k)$ for which 
$\k$ depends continuously on some complex parameters. 
All of them are simple solutions of degree $d \le 4$. 
Except for these solutions, it seems that assertion (1) 
remains true for all non-Riccati algebraic solutions of degree $d \le 6$, 
although a further check is needed to swear its truth (see 
also Remark \ref{rem:low-degree}). 
\item Assertion (2) is not necessarily true if the solution has 
a univalent local branch at a fixed singular point. 
This can be seen by the ``generic" icosahedral solution 
of Boalch \cite{Boalch2}, 
\[
\left\{
\begin{array}{rcl}
z &=& \dfrac{27 s^5(s^2+1)^2(3s-4)^3}{4(2s-1)^3(9s^2+4)^2}, 
\\[6mm]
q &=& \dfrac{3s(3s-4)(s^2+1)(3s^2-2s+4)}{2(2s-1)^2(9s^2+4)}, 
\\[6mm]
p &=& -\dfrac{(2s-1)^2(9s^2+4)
(9s^2+3s+10)}{90 s(3s-4)(s^2+1)(3s^2-3s+2)(3s^2+2s+2)}, 
\end{array}
\right.
\]
for which $d = 12$ and $\k = (1/5,11/60,17/60,7/60,1/60)$. 
This solution has ramification indices (partitions of $d = 12$): 
$5$, $3$, $2$, $2$ at $z = 0$, $\infty$; and 
$3$, $3$, $2$, $2$, $1$, $1$ at $z = 1$. 
So it has two univalent local branches at $z = 1$. 
Observe that $d \k_i$ $(i = 0,1,2,3,4)$ are {\sl not} integers. 
We remark that assertion (2) is valid for any $d \ge 1$ 
(not only for $d \ge 7$). 
\end{enumerate}
\end{remark}
\section{Trigonometric Diophantine Conditions} \label{sec:tdc}
The rationality result in \S\ref{sec:rationality} is stated in 
the Painlev\'e-equation side. 
Switching to the character-variety side, we present another 
result showing that the coordinates of any finite orbit of 
degree $d \ge 7$ are tied down by very tight conditions, 
namely, by certain trigonometric Diophantine conditions. 
In this section we work on $\Sol(\th)$ downstairs rather than 
$\wt{\Sol}(\th)$ upstairs so that the degree means the 
number of points in the finite $G(2)$-orbit on $\Sol(\th)$ 
under consideration. 
\begin{theorem} \label{thm:diophantine}
Given any $\th = (\th_1,\th_2,\th_3,\th_4) \in \Th$，let 
$\O \subset \Sol(\th)$ be a $(\mbox{possibly infinite})$ 
$G(2)$-orbit of degree $d \ge 7$．
Then the orbit $\O$ is finite if and only if 
\begin{equation} \label{eqn:tdc}
\O \subset \Sol(\th) \cap (2 \cos \pi \Q)^3. 
\end{equation}
If this is the case then $\th_1$, $\th_2$, $\th_3$ and 
$\th_4$ must be real cyclotomic integers such that 
\[
-8 < \th_1, \th_2, \th_3 < 8, \qquad -28 < \th_4 < 28.
\]
\end{theorem}
As a corollary, if $\O \subset \Sol(\th)$ is a finite orbit of 
degree $d \ge 7$, then $\th = (\th_1,\th_2,\th_3,\th_4)$ must be 
real and the orbit $\O$ must lie in the real part $\Sol(\th)_{\bR}$ 
of the complex surface $\Sol(\th)$. 
Thus it is also important to investigate the real dynamics on the 
real cubic surface $\Sol(\th)_{\bR}$ with $\th \in \bR^4$. 
\begin{remark} \label{rem:low-degree} 
All finite orbits of degree $d \le 4$ has been classified by 
Cantat and Loray \cite{CL}. 
During the author's visit to Rennes in March 2008, having heard 
of the author's results for $d \ge 7$, F. Loray carried out computer 
experiments to determine all finite orbits of degrees $5$ and $6$. 
These orbits correspond to some algebraic solutions by Theorem 
\ref{thm:alg=finite} and actually it seems that they correspond to 
already known algebraic solutions (a further careful check is needed). 
\end{remark}
\begin{remark} \label{rem:tde} 
It follows from (\ref{eqn:tdc}) in Theorem \ref{thm:diophantine} that 
enumerating all finite orbits of degree $d \ge 7$ can be embedded 
into the problem of solving the trigonometric Diophantine equation 
\begin{equation} \label{eqn:tde}
\sum_{k=1}^8 \cos \pi \xi_k = 0, \qquad 
\xi = (\xi_1,\dots,\xi_8) \in \Q^8. 
\end{equation}
Similar but more tractable trigonometric Diophantine equations 
have appeared in many places (see e.g. \cite{PR,PX} and the 
references therein). 
Although getting harder, equation (\ref{eqn:tde}) still seems to 
be a tractable problem in computer-assisted mathematics. 
However, even if one succeeds in enumerating all solutions to 
equation (\ref{eqn:tde}), there remains the extra job of 
identifying which solutions are relevant to our original problem. 
In any case the author prefers more insightful geometric 
approaches. 
\end{remark}
\par
The proof of Theorem \ref{thm:diophantine} relies largely on 
the direct manipulations of the dynamics on the character variety, 
but it also depends heavily on Theorem \ref{thm:rationality}, which 
in turn is obtained by the combination of some main discussions on 
the Painlev\'e-equation side and some auxiliary discussions on 
the character-variety side. 
Behind this complicated circle of ideas, there exists the 
geometry of cubic surfaces, especially the configuration 
of lines on a cubic surface. 
In the next section we give a brief account of this, 
leaving a full explanation in \cite{Iwasaki2}. 
\section{Lines on a Cubic Surface} \label{sec:lines} 
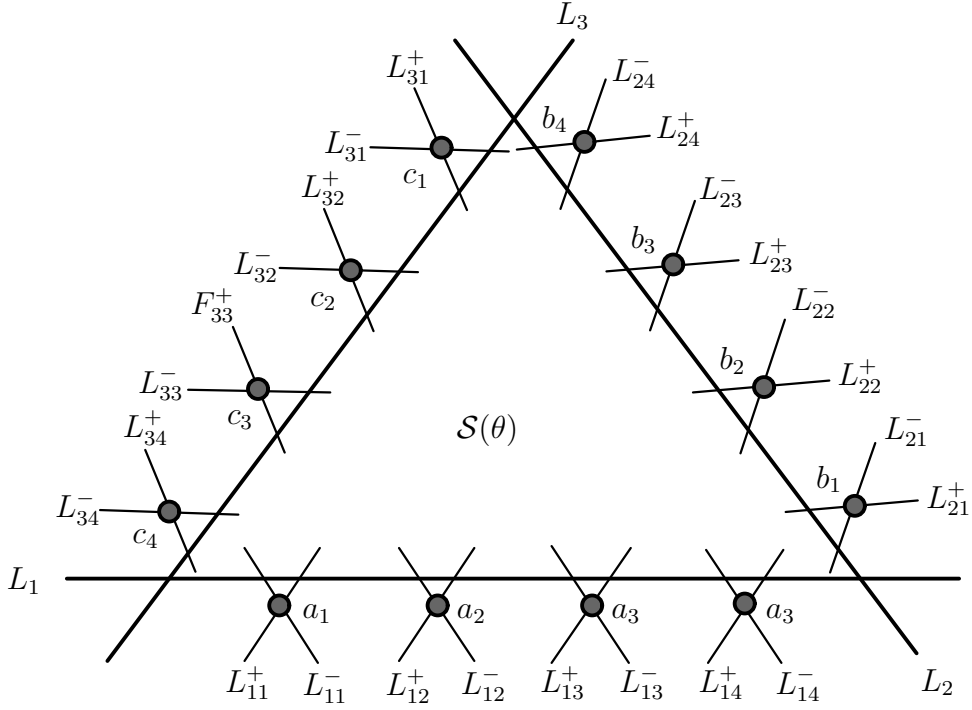
\begin{figure}[t]
\begin{center}
\unitlength 0.1in
\begin{picture}( 49.3000, 35.2000)(  8.6000,-39.6000)
\put(8.6000,-35.6000){\makebox(0,0)[lb]{$L_1$}}%
\put(55.9000,-40.9000){\makebox(0,0)[lb]{$L_2$}}%
\put(37.0000,-6.1000){\makebox(0,0)[lb]{$L_3$}}%
%
\special{pn 20}%
\special{pa 1170 3480}%
\special{pa 5790 3480}%
\special{fp}%
%
\special{pn 20}%
\special{pa 3790 670}%
\special{pa 1380 3910}%
\special{fp}%
\special{pa 3180 670}%
\special{pa 5570 3870}%
\special{fp}%
%
\special{pn 13}%
\special{pa 2090 3320}%
\special{pa 2470 3920}%
\special{fp}%
\special{pa 2480 3320}%
\special{pa 2090 3910}%
\special{fp}%
\special{pa 2890 3320}%
\special{pa 3280 3910}%
\special{fp}%
\special{pa 3280 3320}%
\special{pa 2890 3910}%
\special{fp}%
\special{pa 3680 3310}%
\special{pa 4080 3920}%
\special{fp}%
\special{pa 4090 3320}%
\special{pa 3680 3920}%
\special{fp}%
\special{pa 4480 3330}%
\special{pa 4870 3920}%
\special{fp}%
\special{pa 4880 3330}%
\special{pa 4490 3920}%
\special{fp}%
%
\special{pn 13}%
\special{pa 3500 1240}%
\special{pa 4188 1168}%
\special{fp}%
\special{pa 3726 1550}%
\special{pa 3960 874}%
\special{fp}%
\special{pa 3964 1876}%
\special{pa 4648 1818}%
\special{fp}%
\special{pa 4190 2184}%
\special{pa 4422 1508}%
\special{fp}%
\special{pa 4412 2508}%
\special{pa 5118 2446}%
\special{fp}%
\special{pa 4658 2826}%
\special{pa 4888 2128}%
\special{fp}%
\special{pa 4892 3132}%
\special{pa 5576 3072}%
\special{fp}%
\special{pa 5122 3448}%
\special{pa 5356 2772}%
\special{fp}%
%
\special{pn 13}%
\special{pa 2970 920}%
\special{pa 3242 1556}%
\special{fp}%
\special{pa 2742 1228}%
\special{pa 3458 1250}%
\special{fp}%
\special{pa 2502 1550}%
\special{pa 2758 2188}%
\special{fp}%
\special{pa 2272 1858}%
\special{pa 2988 1880}%
\special{fp}%
\special{pa 2030 2166}%
\special{pa 2298 2822}%
\special{fp}%
\special{pa 1798 2494}%
\special{pa 2534 2508}%
\special{fp}%
\special{pa 1576 2808}%
\special{pa 1836 3444}%
\special{fp}%
\special{pa 1342 3122}%
\special{pa 2058 3146}%
\special{fp}%
\put(56.2000,-31.5000){\makebox(0,0)[lb]{$L_{21}^{+}$}}%
\put(54.0000,-27.9000){\makebox(0,0)[lb]{$L_{21}^{-}$}}%
\put(51.6000,-25.1000){\makebox(0,0)[lb]{$L_{22}^{+}$}}%
\put(49.2000,-21.2000){\makebox(0,0)[lb]{$L_{22}^{-}$}}%
\put(47.0000,-18.8000){\makebox(0,0)[lb]{$L_{23}^{+}$}}%
\put(44.4000,-15.4000){\makebox(0,0)[lb]{$L_{23}^{-}$}}%
\put(42.2000,-12.2000){\makebox(0,0)[lb]{$L_{24}^{+}$}}%
\put(39.9000,-9.2000){\makebox(0,0)[lb]{$L_{24}^{-}$}}%
\put(28.2000,-9.0000){\makebox(0,0)[lb]{$L_{31}^{+}$}}%
\put(25.0000,-13.0000){\makebox(0,0)[lb]{$L_{31}^{-}$}}%
\put(23.8000,-15.3000){\makebox(0,0)[lb]{$L_{32}^{+}$}}%
\put(20.3000,-19.3000){\makebox(0,0)[lb]{$L_{32}^{-}$}}%
\put(18.1000,-21.4000){\makebox(0,0)[lb]{$F_{33}^{+}$}}%
\put(15.4000,-25.5000){\makebox(0,0)[lb]{$L_{33}^{-}$}}%
\put(14.6000,-27.8000){\makebox(0,0)[lb]{$L_{34}^{+}$}}%
\put(11.1000,-32.0000){\makebox(0,0)[lb]{$L_{34}^{-}$}}%
\put(19.9000,-41.1000){\makebox(0,0)[lb]{$L_{11}^+$}}%
\put(23.9000,-41.3000){\makebox(0,0)[lb]{$L_{11}^-$}}%
\put(28.2000,-41.3000){\makebox(0,0)[lb]{$L_{12}^+$}}%
\put(32.1000,-41.2000){\makebox(0,0)[lb]{$L_{12}^-$}}%
\put(36.2000,-41.1000){\makebox(0,0)[lb]{$L_{13}^+$}}%
\put(40.3000,-41.1000){\makebox(0,0)[lb]{$L_{13}^-$}}%
\put(44.4000,-41.1000){\makebox(0,0)[lb]{$L_{14}^+$}}%
\put(48.4000,-41.2000){\makebox(0,0)[lb]{$L_{14}^-$}}%
\put(31.9000,-27.9000){\makebox(0,0)[lb]{$\Sol(\th)$}}%
%
\special{pn 20}%
\special{sh 0.600}%
\special{ar 3850 1200 52 52  0.0000000 6.2831853}%
%
\special{pn 20}%
\special{sh 0.600}%
\special{ar 4310 1840 52 52  0.0000000 6.2831853}%
%
\special{pn 20}%
\special{sh 0.600}%
\special{ar 4780 2480 52 52  0.0000000 6.2831853}%
%
\special{pn 20}%
\special{sh 0.600}%
\special{ar 5250 3100 52 52  0.0000000 6.2831853}%
%
\special{pn 20}%
\special{sh 0.600}%
\special{ar 2270 3620 52 52  0.0000000 6.2831853}%
%
\special{pn 20}%
\special{sh 0.600}%
\special{ar 3090 3620 52 52  0.0000000 6.2831853}%
%
\special{pn 20}%
\special{sh 0.600}%
\special{ar 3890 3620 52 52  0.0000000 6.2831853}%
%
\special{pn 20}%
\special{sh 0.600}%
\special{ar 4680 3610 52 52  0.0000000 6.2831853}%
%
\special{pn 20}%
\special{sh 0.600}%
\special{ar 3110 1230 52 52  0.0000000 6.2831853}%
%
\special{pn 20}%
\special{sh 0.600}%
\special{ar 2640 1870 52 52  0.0000000 6.2831853}%
%
\special{pn 20}%
\special{sh 0.600}%
\special{ar 2160 2490 52 52  0.0000000 6.2831853}%
%
\special{pn 20}%
\special{sh 0.600}%
\special{ar 1700 3130 52 52  0.0000000 6.2831853}%
\put(23.9000,-37.1000){\makebox(0,0)[lb]{$a_1$}}%
\put(31.9000,-37.1000){\makebox(0,0)[lb]{$a_2$}}%
\put(39.9000,-37.1000){\makebox(0,0)[lb]{$a_3$}}%
\put(47.9000,-37.1000){\makebox(0,0)[lb]{$a_3$}}%
\put(50.5000,-30.3000){\makebox(0,0)[lb]{$b_1$}}%
\put(45.5000,-24.1000){\makebox(0,0)[lb]{$b_2$}}%
\put(40.9000,-17.8000){\makebox(0,0)[lb]{$b_3$}}%
\put(36.3000,-11.5000){\makebox(0,0)[lb]{$b_4$}}%
\put(15.1000,-33.2000){\makebox(0,0)[lb]{$c_4$}}%
\put(29.1000,-14.4000){\makebox(0,0)[lb]{$c_1$}}%
\put(24.3000,-20.8000){\makebox(0,0)[lb]{$c_2$}}%
\put(20.0000,-27.0000){\makebox(0,0)[lb]{$c_3$}}%
\end{picture}%
\end{center}
\caption{The 27 lines viewed from the tritangent lines at 
infinity}
\label{fig:lines}
\end{figure} 
Compactify the affine cubic surface $\Sol(\th)$ by 
the standard embedding $\Sol(\th) \hookrightarrow 
\ol{\Sol}(\th) \subset \P^3$. 
Then $\ol{\Sol}(\th)$ is obtained from $\Sol(\th)$ 
by adding the tritangent lines at infinity, 
$L = L_1 \cup L_2 \cup L_3$, as in Figure \ref{fig:lines}. 
For simplicity we assume that $\th = \rh(\k)$ with 
$\k \in \K - \Wall(D_4^{(1)})$. 
Then the projective cubic surface $\ol{\Sol}(\th)$ is 
smooth and it contains twenty-seven lines, whose configuration 
is depicted in Figure \ref{fig:lines}. 
The lines at infinity, $L_1$, $L_2$, $L_3$, are three among them. 
The remaining twenty-four lines are divided into three groups, each 
consisting of eight lines, according to the three lines at infinity. 
Namely, for each $i = 1,2,3$, the line $L_i$ meets exactly 
eight lines, say, $L_{ij}^{\ve}$ as in Figure \ref{fig:lines}, 
where $j = 1,2,3,4$ and $\ve = \pm$. 
This group of eight lines are divided into four intersecting 
pairs $\{L_{ij}^+, L_{ij}^{-}\}_{j=1}^4$. 
Any other pair from the same group has no intersections. 
\par
Assume that a finite $G(2)$-orbit $\O \subset \Sol(\th)$ 
be given. 
To it we can associate an ``ON/OFF" data 
$(\ab,\bb,\cb) \in \{0,1 \}^{12}$ as follows. 
To define $\ab = (a_1,a_2,a_3,a_4) \in \{0,1\}^4$, we put 
\[
a_j := \left\{
\begin{array}{ll}
1 \,\, \mbox{(ON)}, \quad & \mbox{if $\O$ passes through the 
intersection point $L_{1j}^{+} \cap L_{1j}^{-}$}, 
\\[2mm] 
0 \,\, \mbox{(OFF)}, \quad & \mbox{otherwise}, 
\end{array}
\right.
\]
for $j = 1,2,3,4$. 
In a similar manner we can define 
$\bb = (b_1,b_2,b_3,b_4) \in \{0,1\}^4$ and 
$\cb = (c_1,c_2,c_3,c_4) \in \{0,1\}^4$ by replacing 
$L_{1j}^{\pm}$ with $L_{2j}^{\pm}$ and $L_{3j}^{\pm}$ 
repectively. 
Then certain arguments that are too involved to be included here 
lead to the matrix 
\[
M(\ab,\bb,\cb) := 
\begin{pmatrix} 
d_1 & a_3-a_4 & c_1-c_2 & b_1-b_2 \\
a_3-a_4 & d_2 & b_3-b_4 & c_3-c_4 \\
c_1-c_2 & b_3-b_4 & d_3 & a_1-a_2 \\
b_1-b_2 & c_3-c_4 & a_1-a_2 & d_4 
\end{pmatrix}, 
\]
where $d_i$ $(i = 1,2,3,4)$ are nonnegative integers 
defined by 
\[
\left\{
\begin{array}{rcl}
d_1 &:=& a_3+a_4+b_1+b_2+c_1+c_2, \\
d_2 &:=& a_3+a_4+b_3+b_4+c_3+c_4, \\
d_3 &:=& a_1+a_2+b_3+b_4+c_1+c_2, \\
d_4 &:=& a_1+a_2+b_1+b_2+c_3+c_4. 
\end{array}
\right.
\]
It turns out that the column vector 
$\kb = {}^t(\k_1,\k_2,\k_3,\k_4)$ must satisfy 
a linear equation 
\begin{equation} \label{eqn:linear} 
[d I_4 - M(\ab,\bb,\cb)] \, \kb = \mbox{a certain integer vector}, 
\end{equation}
where $d$ is the order of the orbit $\O$ and $I_4$ is the 
identity matrix of rank $4$. 
\begin{lemma} \label{lem:matrix} 
For any $(\ab,\bb,\cb) \in \{0,1\}^{12}$ the matrix 
$M(\ab,\bb,\cb)$ has no eigenvalues $\ge 7$. 
\end{lemma}
This is verified by a computer check exhausting all $2^{12} = 4096$ 
possibilities for the data $(\ab,\bb,\cb)$. 
It is also observed that actually some of $0,1,2,3,4,5,6$ are  
eigenvalues of the matrix $M(\ab,\bb,\cb)$. 
The author is indebted to A. Maruyama and T. Uehara for the job of 
these verifications. 
\par\medskip\noindent 
{\bf Sketch of the proof of Theorem \ref{thm:rationality}}. 
If $d \ge 7$ then Lemma \ref{lem:matrix} implies that 
$d I_4 - M(\ab,\bb,\cb)$ is invertible in rational numbers 
since it is an integer matrix, so that equation (\ref{eqn:linear}) 
can be settled to conclude that $\kb$ is a vector with rational 
entries. This proves assertion (1). 
Let us proceed to assertion (2). 
Put $z_1 = 0$, $z_2 = 1$ and $z_3 = \infty$. 
It is shown in \cite{Iwasaki1} that for each $i = 1,2,3$ the line 
$L_i$ at infinity is attached to the fixed singular point $z_i$ and 
the univalent solution germs at $z_i$ are in one-to-one 
correspondence with those intersection points 
$L_{ij}^{+} \cap L_{ij}^{-}$, $j \in \{1,2,3,4\}$, which 
lie in the affine part $\Sol(\th)$ of $\ol{\Sol}(\th)$. 
Thus if the algebraic solution under consideration has no 
univalent local branches at any fixed singular point, then 
we must have $(\ab,\bb,\cb)=(\0,\0,\0)$ and 
$M(\ab,\bb,\cb) = O$. 
Then equation (\ref{eqn:linear}) implies that 
$d \kb$ must be an integer vector, from which assertion (2) 
readily follows. 
Note that this argument is valid for an arbitrary integer 
$d \ge 1$. \hfill $\Box$
\par\medskip
This section ends with three remarks. 
Firstly, even if $d \le 6$ some useful information about 
$\kb$ can be extracted from equation (\ref{eqn:linear}). 
Secondly, if $\k \in \Wall(D_4^{(1)})$ then the line 
configuration is degenerate and the situation becomes more 
complicated than the case $\k \in \K-\Wall(D_4^{(1)})$ 
discussed above, but basically a similar argument is feasible. 
Finally we refer to the original paper \cite{Iwasaki2} 
for the most important thing: why and how equation 
(\ref{eqn:linear}) occurs. 
\section{Stratifications of Parameters} \label{sec:strata}
\begin{figure}[t] 
\begin{center}
\unitlength 0.1in
\begin{picture}( 51.7000, 15.1000)(  2.7000,-21.7000)
%
\special{pn 20}%
\special{sh 0.600}%
\special{ar 448 852 48 48  0.0000000 6.2831853}%
%
\special{pn 20}%
\special{sh 0.600}%
\special{ar 1248 844 48 48  0.0000000 6.2831853}%
%
\special{pn 20}%
\special{sh 0.600}%
\special{ar 456 1652 48 48  0.0000000 6.2831853}%
%
\special{pn 20}%
\special{ar 1256 1652 48 48  0.0000000 6.2831853}%
%
\special{pn 20}%
\special{sh 0.600}%
\special{ar 848 1252 48 48  0.0000000 6.2831853}%
%
\special{pn 20}%
\special{pa 480 900}%
\special{pa 816 1228}%
\special{fp}%
\special{pa 824 1228}%
\special{pa 824 1228}%
\special{fp}%
%
\special{pn 20}%
\special{pa 880 1292}%
\special{pa 1216 1620}%
\special{dt 0.054}%
\special{pa 1216 1620}%
\special{pa 1216 1620}%
\special{dt 0.054}%
%
\special{pn 20}%
\special{pa 1216 884}%
\special{pa 880 1212}%
\special{fp}%
%
\special{pn 20}%
\special{pa 808 1292}%
\special{pa 504 1604}%
\special{fp}%
\put(2.7000,-8.3000){\makebox(0,0)[lb]{$1$}}%
\put(13.5000,-8.4000){\makebox(0,0)[lb]{$2$}}%
\put(2.7000,-18.2000){\makebox(0,0)[lb]{$3$}}%
\put(13.2800,-18.4400){\makebox(0,0)[lb]{$4$}}%
\put(7.9000,-11.5000){\makebox(0,0)[lb]{$0$}}%
%
\special{pn 20}%
\special{sh 0.600}%
\special{ar 2536 852 48 48  0.0000000 6.2831853}%
%
\special{pn 20}%
\special{sh 0.600}%
\special{ar 3336 844 48 48  0.0000000 6.2831853}%
%
\special{pn 20}%
\special{sh 0.600}%
\special{ar 2544 1652 48 48  0.0000000 6.2831853}%
%
\special{pn 20}%
\special{sh 0.600}%
\special{ar 3344 1652 48 48  0.0000000 6.2831853}%
%
\special{pn 20}%
\special{ar 2936 1252 48 48  0.0000000 6.2831853}%
%
\special{pn 20}%
\special{pa 2568 900}%
\special{pa 2904 1228}%
\special{dt 0.054}%
\special{pa 2912 1228}%
\special{pa 2912 1228}%
\special{dt 0.054}%
%
\special{pn 20}%
\special{pa 2968 1292}%
\special{pa 3304 1620}%
\special{dt 0.054}%
\special{pa 3304 1620}%
\special{pa 3304 1620}%
\special{dt 0.054}%
%
\special{pn 20}%
\special{pa 3304 884}%
\special{pa 2968 1212}%
\special{dt 0.054}%
%
\special{pn 20}%
\special{pa 2896 1292}%
\special{pa 2592 1604}%
\special{dt 0.054}%
\put(23.7000,-8.3000){\makebox(0,0)[lb]{$1$}}%
\put(34.5000,-8.4000){\makebox(0,0)[lb]{$2$}}%
\put(23.7000,-17.9000){\makebox(0,0)[lb]{$3$}}%
\put(34.1600,-18.3600){\makebox(0,0)[lb]{$4$}}%
\put(28.8000,-11.5000){\makebox(0,0)[lb]{$0$}}%
%
\special{pn 20}%
\special{sh 0.600}%
\special{ar 4552 850 48 48  0.0000000 6.2831853}%
%
\special{pn 20}%
\special{sh 0.600}%
\special{ar 5352 842 48 48  0.0000000 6.2831853}%
%
\special{pn 20}%
\special{ar 4560 1650 48 48  0.0000000 6.2831853}%
%
\special{pn 20}%
\special{ar 5360 1650 48 48  0.0000000 6.2831853}%
%
\special{pn 20}%
\special{sh 0.600}%
\special{ar 4952 1250 48 48  0.0000000 6.2831853}%
%
\special{pn 20}%
\special{pa 4584 898}%
\special{pa 4920 1226}%
\special{fp}%
\special{pa 4928 1226}%
\special{pa 4928 1226}%
\special{fp}%
%
\special{pn 20}%
\special{pa 4984 1290}%
\special{pa 5320 1618}%
\special{dt 0.054}%
\special{pa 5320 1618}%
\special{pa 5320 1618}%
\special{dt 0.054}%
%
\special{pn 20}%
\special{pa 5320 882}%
\special{pa 4984 1210}%
\special{fp}%
%
\special{pn 20}%
\special{pa 4912 1290}%
\special{pa 4608 1602}%
\special{dt 0.054}%
\put(43.6000,-8.3000){\makebox(0,0)[lb]{$1$}}%
\put(54.4000,-8.3000){\makebox(0,0)[lb]{$2$}}%
\put(43.8000,-18.2000){\makebox(0,0)[lb]{$3$}}%
\put(54.0800,-18.4200){\makebox(0,0)[lb]{$4$}}%
\put(49.0000,-11.5000){\makebox(0,0)[lb]{$0$}}%
\put(7.0400,-20.0400){\makebox(0,0)[lb]{$D_4$}}%
\put(27.8400,-20.5200){\makebox(0,0)[lb]{$A_1^{\oplus 4}$}}%
\put(48.3200,-20.1000){\makebox(0,0)[lb]{$A_3$}}%
\put(2.9600,-23.3200){\makebox(0,0)[lb]{$I = \{0,1,2,3\}$}}%
\put(24.0800,-23.3200){\makebox(0,0)[lb]{$I = \{1,2,3,4\}$}}%
\put(44.7000,-23.4000){\makebox(0,0)[lb]{$I = \{0,1,2\}$}}%
\end{picture}%
\end{center}
\caption{Some $D_4^{(1)}$-strata and their abstract Dynkin types}
\label{fig:strata}
\end{figure}
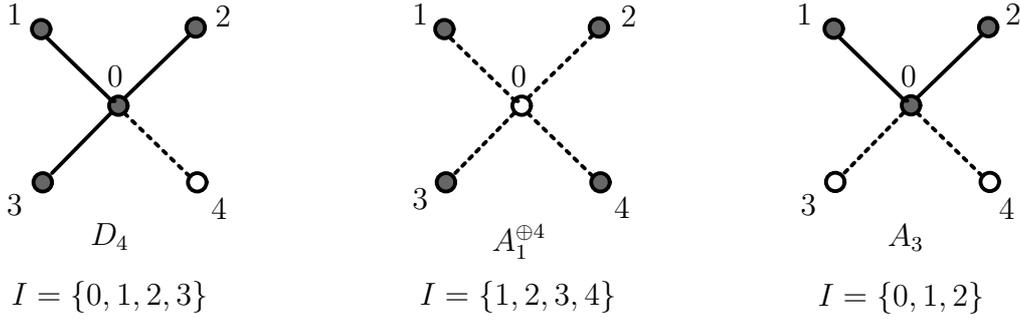

We define a stratification of the parameter space $\K$ in terms 
the proper subdiagrams of the Dynkin diagram $D_4^{(1)}$. 
To this end we index the nodes of the Dynkin diagram $D_4^{(1)}$ 
by the numbers $0, 1, 2, 3, 4$, where $0$ represents 
the central node (see Figure \ref{fig:strata}). 
Let $\I$ be the set of all proper subsets of $\{0,1,2,3,4\}$ 
including the empty set $\emptyset$. 
For each $I \in \I$ we put 
\begin{equation} \label{eqn:stratif} 
\begin{array}{rcl}
\ol{\K}_I &=& \mbox{the $W(D_4^{(1)})$-translates of the subset} \,\, 
\{\, \k \in \K \,:\, \k_i = 0 \,\, (i \in I) \,\}, 
\\[2mm]
\K_I &=& \ol{\K}_I - \displaystyle \bigcup_{|J|=|I|+1} \ol{\K}_J, 
\\[2mm]
D_I &=& \mbox{the Dynkin subdiagram of $D_4^{(1)}$ that has nodes 
$\bullet$ exactly in $I$}. 
\end{array}
\end{equation}
It turns out that for any pair $(I, I') \in \I \times \I$, either 
$\K_I = \K_{I'}$ or $\K_I \cap \K_{I'} = \emptyset$ holds so that 
the partition $\{\K_I\}_{I \in \I}$ defines a stratification of 
$\K$, called the $D_4^{(1)}$-{\sl stratification}. 
For $I = \emptyset$ one has the big open stratum 
$\K_{\emptyset} = \K - \Wall(D_4^{(1)})$ and other examples 
of strata are given in Figure \ref{fig:strata}. 
\par
The automorphism group of Dynkin diagram $D_4^{(1)}$ is 
the symmetric group $S_4$ of degree $4$ acting by permuting 
the nodes $1$, $2$, $3$, $4$ while fixing the central node $0$. 
The group $W(D_4^{(1)})$ extended by $S_4$ is the affine 
Weyl group $W(F_4^{(1)})$ of type $F_4^{(1)}$. 
A coaser stratification of $\K$ can be defined in the same way 
as in the case of $D_4^{(1)}$-stratification  
by replacing the group $W(D_4^{(1)})$ with $W(F_4^{(1)})$ 
in (\ref{eqn:stratif}). 
It is called the $F_4^{(1)}$-{\sl stratification}. 
Note that the $F_4^{(1)}$-stratification encodes only the 
{\sl abstract} Dynkin type of the subdiagram $D_I$, while the 
$D_4^{(1)}$-stratification encodes not only the abstract Dynkin type of 
$D_{I}$ but also the inclusion patern $D_I \hookrightarrow D_4^{(1)}$, 
a kind of marking. 
Thus the $F_4^{(1)}$-strata can be indexed by the abstract Dynkin 
subdiagrams of $D_4^{(1)}$. 
The adjacency relations among them are given in Figure \ref{fig:adjac}, 
where $* \to **$ indicates that the stratum $**$ is in the closure 
of $*$. 
\begin{figure}[t]
\[
\begin{CD}
\emptyset @>>> A_1 @>>> A_1^{\oplus 2} @>>> A_1^{\oplus 3} 
@>>> A_1^{\oplus 4} \\
@. @VVV @VVV @VVV @. \\
 @. A_2 @>>> A_3 @>>> D_4 @.  
\end{CD}
\]
\caption{Adjacency relations among the $F_4^{(1)}$-strata}
\label{fig:adjac}
\end{figure}
\section{On Various Strata} \label{sec:various} 
Theorems \ref{thm:rationality} and \ref{thm:diophantine} are results 
that can be stated without refering to the stratification. 
Besides them, there are such results that differ stratum by stratum. 
A factor that might cause such a difference is the topology 
(or perhaps the shape) of the real character variety $\Sol(\th)_{\bR}$ 
(see \cite{BeGo}). 
On one hand the topology changes as the stratum varies and 
on the other hand the dynamics of the mapping class group action on 
$\Sol(\th)_{\bR}$ is {\sl a priori} defined by the space 
$\Sol(\th)_{\bR}$ itself, so that the topology or the shape of 
the space should have a strong influence on the dynamics. 
\par
\begin{figure}[h]
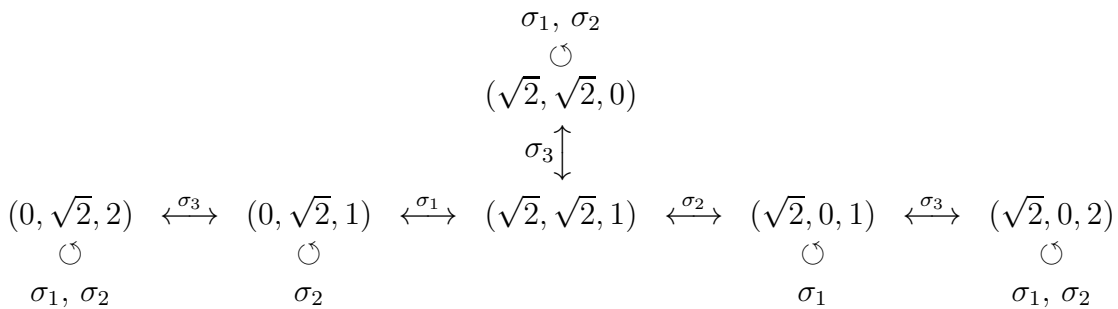

\[
\begin{array}{ccccccccc}
 & & & & \si_1, \, \si_2 & & & & \\
 & & & & \carl & & & & \\
 & & & & (\sqrt{2},\sqrt{2},0) & & & & \\[1mm]
 & & & & 
\mbox{\rot{90}{$\stackrel{\mbox{\rot{-90}{$\si_3$}}}{\llra}$}} 
\phantom{\si_3} & & & & \\
(0,\sqrt{2},2) & \stackrel{\si_3}{\llra} & 
(0,\sqrt{2},1) & \stackrel{\si_1}{\llra} & 
(\sqrt{2},\sqrt{2},1) & \stackrel{\si_2}{\llra} & 
(\sqrt{2},0,1) & \stackrel{\si_3}{\llra} & 
(\sqrt{2},0,2) \\
\carl & & \carl & & & & \carl & & \carl \\
\si_1, \, \si_2 & & \si_2 & & & & \si_1 & & \si_1, \, \si_2 
\end{array}
\]
\caption{A finite orbit of degree $6$ on the stratum of type 
$A_1^{\oplus2}$}
\label{fig:A12}
\end{figure}

We focus our attention on the $F_4^{(1)}$-strata of positive 
codimensions. 
A careful inspection shows that it is natural to divide those strata 
into two sequences (see Figure \ref{fig:adjac}): 
\begin{center}
(S1) \quad $A_1^{\oplus2} \to A_1^{\oplus3}\to A_1^{\oplus4}$, 
\qquad 
(S2) \quad $A_1 \to A_2 \to A_3 \to D_4$. 
\end{center}
In this section we are concerned with the strata belonging to 
the former sequence (S1). 
\begin{example}[Stratum of type $\mbox{\boldmath $A_1^{\oplus4}$}$] 
\label{ex:A14}
This is the locus where the classically well-known Picard solutions 
exist (see \cite{Mazzocco2})．
The corresponding character variety $\Sol(\th)$ is the Cayley cubic, 
with parameters $\th = (0,0,0,-4)$. 
The Picard solutions can be settled by quadrature in terms of the 
Legendre family of elliptic curves. 
However the way in which they are integrated is irreducible in the 
sense of Nishioka \cite{Nishioka} and Umemura \cite{Umemura}, but 
reducible in the sense of Casale \cite{Casale} and 
Malgrange \cite{Malgrange} (see Cantat and Loray \cite{CL})．
This world is amenable to torus structures in two ways. 
Firstly an elliptic curve is a (real) torus and secondly the Cayley 
cubic enjoys a (complex) orbifold torus structure，
$\Sol(\th) \cong (\C^{\times})^2/(\mbox{an involution})$, where the 
four $A_1$-singularities (all real) just come from the four fixed points 
of the involution. 
On this stratum there are countably many algebraic solutions, 
which correspond to the finite-order points of elliptic curves．
The finite orbits on the Cayley cubic are dense in the unique 
bounded connected component of the real Cayley cubic 
$\Sol(\th)_{\bR}$ with the four singular points removed 
(see Figure \ref{fig:cayley}). 
\end{example}
\begin{figure}[t]
\begin{center}
\unitlength 0.1in
\begin{picture}( 30.5600, 22.6000)(  9.6400,-24.0100)
%
\special{pn 20}%
\special{pa 1970 230}%
\special{pa 2000 218}%
\special{pa 2032 208}%
\special{pa 2062 196}%
\special{pa 2092 186}%
\special{pa 2122 176}%
\special{pa 2154 168}%
\special{pa 2184 160}%
\special{pa 2214 154}%
\special{pa 2246 148}%
\special{pa 2278 144}%
\special{pa 2308 142}%
\special{pa 2340 142}%
\special{pa 2372 142}%
\special{pa 2404 142}%
\special{pa 2434 144}%
\special{pa 2466 148}%
\special{pa 2498 152}%
\special{pa 2530 156}%
\special{pa 2562 162}%
\special{pa 2594 168}%
\special{pa 2628 174}%
\special{pa 2660 180}%
\special{pa 2692 188}%
\special{pa 2724 194}%
\special{pa 2750 200}%
\special{sp}%
%
\special{pn 20}%
\special{pa 1970 240}%
\special{pa 2002 250}%
\special{pa 2034 260}%
\special{pa 2066 270}%
\special{pa 2096 278}%
\special{pa 2128 288}%
\special{pa 2160 296}%
\special{pa 2192 302}%
\special{pa 2222 308}%
\special{pa 2254 314}%
\special{pa 2286 318}%
\special{pa 2316 320}%
\special{pa 2348 322}%
\special{pa 2380 322}%
\special{pa 2410 320}%
\special{pa 2442 316}%
\special{pa 2472 310}%
\special{pa 2504 304}%
\special{pa 2534 296}%
\special{pa 2564 288}%
\special{pa 2596 278}%
\special{pa 2626 268}%
\special{pa 2658 256}%
\special{pa 2688 244}%
\special{pa 2718 232}%
\special{pa 2750 220}%
\special{pa 2750 220}%
\special{sp}%
%
\special{pn 20}%
\special{pa 2750 230}%
\special{pa 2734 260}%
\special{pa 2718 288}%
\special{pa 2702 318}%
\special{pa 2686 346}%
\special{pa 2670 374}%
\special{pa 2652 402}%
\special{pa 2634 428}%
\special{pa 2616 454}%
\special{pa 2596 480}%
\special{pa 2576 504}%
\special{pa 2556 528}%
\special{pa 2534 550}%
\special{pa 2512 572}%
\special{pa 2486 592}%
\special{pa 2462 610}%
\special{pa 2436 628}%
\special{pa 2408 644}%
\special{pa 2380 660}%
\special{pa 2352 676}%
\special{pa 2324 690}%
\special{pa 2294 704}%
\special{pa 2266 720}%
\special{pa 2236 734}%
\special{pa 2208 748}%
\special{pa 2180 764}%
\special{pa 2152 780}%
\special{pa 2124 798}%
\special{pa 2098 816}%
\special{pa 2072 834}%
\special{pa 2048 854}%
\special{pa 2024 876}%
\special{pa 2002 898}%
\special{pa 1980 920}%
\special{pa 1960 944}%
\special{pa 1938 970}%
\special{pa 1918 994}%
\special{pa 1898 1020}%
\special{pa 1880 1046}%
\special{pa 1860 1072}%
\special{pa 1842 1100}%
\special{pa 1824 1126}%
\special{pa 1806 1154}%
\special{pa 1800 1160}%
\special{sp}%
%
\special{pn 20}%
\special{pa 1980 240}%
\special{pa 1998 270}%
\special{pa 2014 298}%
\special{pa 2030 328}%
\special{pa 2048 356}%
\special{pa 2064 384}%
\special{pa 2082 412}%
\special{pa 2100 438}%
\special{pa 2118 466}%
\special{pa 2138 492}%
\special{pa 2158 516}%
\special{pa 2178 540}%
\special{pa 2200 564}%
\special{pa 2222 586}%
\special{pa 2244 608}%
\special{pa 2268 628}%
\special{pa 2294 646}%
\special{pa 2320 664}%
\special{pa 2348 680}%
\special{pa 2376 694}%
\special{pa 2404 708}%
\special{pa 2434 722}%
\special{pa 2464 736}%
\special{pa 2494 748}%
\special{pa 2524 762}%
\special{pa 2554 774}%
\special{pa 2584 788}%
\special{pa 2612 802}%
\special{pa 2640 818}%
\special{pa 2668 834}%
\special{pa 2694 850}%
\special{pa 2718 870}%
\special{pa 2742 890}%
\special{pa 2764 912}%
\special{pa 2786 934}%
\special{pa 2806 958}%
\special{pa 2824 984}%
\special{pa 2842 1010}%
\special{pa 2860 1038}%
\special{pa 2876 1066}%
\special{pa 2892 1094}%
\special{pa 2908 1124}%
\special{pa 2922 1154}%
\special{pa 2938 1184}%
\special{pa 2950 1210}%
\special{sp}%
%
\special{pn 20}%
\special{pa 1800 1170}%
\special{pa 1824 1194}%
\special{pa 1846 1218}%
\special{pa 1870 1242}%
\special{pa 1894 1266}%
\special{pa 1916 1288}%
\special{pa 1940 1312}%
\special{pa 1964 1334}%
\special{pa 1988 1354}%
\special{pa 2014 1376}%
\special{pa 2038 1396}%
\special{pa 2064 1414}%
\special{pa 2090 1432}%
\special{pa 2116 1450}%
\special{pa 2144 1464}%
\special{pa 2172 1480}%
\special{pa 2200 1492}%
\special{pa 2230 1504}%
\special{pa 2260 1514}%
\special{pa 2290 1524}%
\special{pa 2322 1532}%
\special{pa 2352 1540}%
\special{pa 2382 1550}%
\special{pa 2414 1560}%
\special{pa 2444 1572}%
\special{pa 2472 1584}%
\special{pa 2502 1600}%
\special{pa 2528 1616}%
\special{pa 2556 1634}%
\special{pa 2582 1654}%
\special{pa 2608 1674}%
\special{pa 2632 1696}%
\special{pa 2656 1718}%
\special{pa 2678 1742}%
\special{pa 2700 1766}%
\special{pa 2720 1792}%
\special{pa 2740 1820}%
\special{pa 2758 1846}%
\special{pa 2776 1874}%
\special{pa 2792 1902}%
\special{pa 2806 1930}%
\special{pa 2820 1960}%
\special{pa 2834 1990}%
\special{pa 2846 2018}%
\special{pa 2858 2048}%
\special{pa 2868 2078}%
\special{pa 2880 2108}%
\special{pa 2890 2138}%
\special{pa 2900 2168}%
\special{pa 2900 2170}%
\special{sp}%
%
\special{pn 20}%
\special{pa 2960 1190}%
\special{pa 2938 1216}%
\special{pa 2916 1240}%
\special{pa 2894 1264}%
\special{pa 2872 1288}%
\special{pa 2850 1310}%
\special{pa 2826 1332}%
\special{pa 2802 1354}%
\special{pa 2778 1374}%
\special{pa 2752 1392}%
\special{pa 2726 1410}%
\special{pa 2700 1426}%
\special{pa 2672 1440}%
\special{pa 2644 1454}%
\special{pa 2616 1466}%
\special{pa 2586 1476}%
\special{pa 2556 1488}%
\special{pa 2524 1498}%
\special{pa 2494 1506}%
\special{pa 2464 1516}%
\special{pa 2432 1526}%
\special{pa 2402 1538}%
\special{pa 2370 1548}%
\special{pa 2340 1560}%
\special{pa 2310 1574}%
\special{pa 2280 1588}%
\special{pa 2250 1604}%
\special{pa 2222 1622}%
\special{pa 2194 1638}%
\special{pa 2166 1658}%
\special{pa 2138 1678}%
\special{pa 2112 1698}%
\special{pa 2088 1720}%
\special{pa 2064 1742}%
\special{pa 2040 1766}%
\special{pa 2018 1790}%
\special{pa 1998 1816}%
\special{pa 1978 1842}%
\special{pa 1960 1868}%
\special{pa 1942 1894}%
\special{pa 1926 1922}%
\special{pa 1912 1952}%
\special{pa 1900 1980}%
\special{pa 1890 2010}%
\special{pa 1880 2042}%
\special{pa 1872 2072}%
\special{pa 1864 2104}%
\special{pa 1858 2134}%
\special{pa 1852 2166}%
\special{pa 1848 2198}%
\special{pa 1842 2230}%
\special{pa 1840 2240}%
\special{sp}%
%
\special{pn 20}%
\special{pa 2940 1200}%
\special{pa 2960 1174}%
\special{pa 2978 1148}%
\special{pa 2998 1122}%
\special{pa 3016 1096}%
\special{pa 3036 1070}%
\special{pa 3058 1046}%
\special{pa 3078 1022}%
\special{pa 3100 1000}%
\special{pa 3124 978}%
\special{pa 3148 958}%
\special{pa 3172 940}%
\special{pa 3198 922}%
\special{pa 3224 906}%
\special{pa 3252 890}%
\special{pa 3280 874}%
\special{pa 3308 860}%
\special{pa 3338 848}%
\special{pa 3368 834}%
\special{pa 3398 822}%
\special{pa 3428 810}%
\special{pa 3460 798}%
\special{pa 3490 788}%
\special{pa 3522 776}%
\special{pa 3554 766}%
\special{pa 3584 756}%
\special{pa 3600 750}%
\special{sp}%
%
\special{pn 20}%
\special{pa 2940 1200}%
\special{pa 2954 1230}%
\special{pa 2968 1260}%
\special{pa 2980 1290}%
\special{pa 2994 1320}%
\special{pa 3010 1350}%
\special{pa 3024 1378}%
\special{pa 3040 1406}%
\special{pa 3056 1432}%
\special{pa 3074 1458}%
\special{pa 3094 1482}%
\special{pa 3114 1506}%
\special{pa 3136 1528}%
\special{pa 3158 1548}%
\special{pa 3182 1568}%
\special{pa 3206 1588}%
\special{pa 3232 1604}%
\special{pa 3258 1622}%
\special{pa 3286 1638}%
\special{pa 3314 1652}%
\special{pa 3342 1668}%
\special{pa 3372 1682}%
\special{pa 3402 1694}%
\special{pa 3432 1708}%
\special{pa 3464 1720}%
\special{pa 3496 1732}%
\special{pa 3528 1744}%
\special{pa 3558 1756}%
\special{pa 3590 1766}%
\special{pa 3624 1778}%
\special{pa 3630 1780}%
\special{sp}%
%
\special{pn 20}%
\special{pa 3590 750}%
\special{pa 3604 780}%
\special{pa 3616 810}%
\special{pa 3628 840}%
\special{pa 3640 868}%
\special{pa 3652 898}%
\special{pa 3664 928}%
\special{pa 3674 958}%
\special{pa 3686 988}%
\special{pa 3696 1020}%
\special{pa 3704 1050}%
\special{pa 3714 1080}%
\special{pa 3722 1112}%
\special{pa 3728 1142}%
\special{pa 3736 1174}%
\special{pa 3740 1206}%
\special{pa 3744 1238}%
\special{pa 3748 1270}%
\special{pa 3750 1304}%
\special{pa 3752 1336}%
\special{pa 3752 1368}%
\special{pa 3750 1400}%
\special{pa 3748 1432}%
\special{pa 3744 1464}%
\special{pa 3740 1496}%
\special{pa 3734 1528}%
\special{pa 3726 1558}%
\special{pa 3718 1588}%
\special{pa 3708 1620}%
\special{pa 3696 1650}%
\special{pa 3684 1678}%
\special{pa 3672 1708}%
\special{pa 3660 1738}%
\special{pa 3646 1768}%
\special{pa 3640 1780}%
\special{sp}%
%
\special{pn 20}%
\special{pa 1840 2230}%
\special{pa 1868 2246}%
\special{pa 1898 2262}%
\special{pa 1926 2278}%
\special{pa 1954 2292}%
\special{pa 1982 2308}%
\special{pa 2012 2322}%
\special{pa 2040 2334}%
\special{pa 2070 2346}%
\special{pa 2100 2356}%
\special{pa 2130 2366}%
\special{pa 2160 2374}%
\special{pa 2190 2382}%
\special{pa 2222 2388}%
\special{pa 2252 2392}%
\special{pa 2284 2396}%
\special{pa 2316 2398}%
\special{pa 2350 2400}%
\special{pa 2382 2400}%
\special{pa 2414 2402}%
\special{pa 2448 2400}%
\special{pa 2482 2400}%
\special{pa 2516 2398}%
\special{pa 2550 2394}%
\special{pa 2584 2390}%
\special{pa 2616 2386}%
\special{pa 2648 2378}%
\special{pa 2680 2370}%
\special{pa 2710 2360}%
\special{pa 2740 2348}%
\special{pa 2768 2332}%
\special{pa 2792 2314}%
\special{pa 2816 2296}%
\special{pa 2838 2272}%
\special{pa 2858 2248}%
\special{pa 2878 2222}%
\special{pa 2896 2194}%
\special{pa 2910 2170}%
\special{sp}%
%
\special{pn 20}%
\special{pa 1800 1180}%
\special{pa 1782 1154}%
\special{pa 1764 1126}%
\special{pa 1746 1100}%
\special{pa 1728 1074}%
\special{pa 1708 1048}%
\special{pa 1688 1022}%
\special{pa 1670 996}%
\special{pa 1648 972}%
\special{pa 1628 948}%
\special{pa 1606 926}%
\special{pa 1584 904}%
\special{pa 1560 882}%
\special{pa 1536 862}%
\special{pa 1512 842}%
\special{pa 1486 824}%
\special{pa 1460 806}%
\special{pa 1432 788}%
\special{pa 1406 772}%
\special{pa 1378 756}%
\special{pa 1350 740}%
\special{pa 1322 724}%
\special{pa 1292 708}%
\special{pa 1264 692}%
\special{pa 1234 678}%
\special{pa 1206 664}%
\special{pa 1180 650}%
\special{sp}%
%
\special{pn 20}%
\special{pa 1800 1190}%
\special{pa 1782 1218}%
\special{pa 1764 1244}%
\special{pa 1746 1272}%
\special{pa 1728 1298}%
\special{pa 1708 1326}%
\special{pa 1688 1350}%
\special{pa 1670 1376}%
\special{pa 1648 1402}%
\special{pa 1628 1424}%
\special{pa 1606 1448}%
\special{pa 1584 1470}%
\special{pa 1562 1492}%
\special{pa 1538 1512}%
\special{pa 1512 1530}%
\special{pa 1486 1548}%
\special{pa 1460 1566}%
\special{pa 1432 1582}%
\special{pa 1404 1596}%
\special{pa 1376 1612}%
\special{pa 1348 1626}%
\special{pa 1318 1638}%
\special{pa 1288 1652}%
\special{pa 1258 1664}%
\special{pa 1228 1676}%
\special{pa 1196 1688}%
\special{pa 1166 1700}%
\special{pa 1140 1710}%
\special{sp}%
%
\special{pn 20}%
\special{pa 1180 650}%
\special{pa 1190 680}%
\special{pa 1198 710}%
\special{pa 1208 738}%
\special{pa 1216 768}%
\special{pa 1226 798}%
\special{pa 1234 828}%
\special{pa 1244 858}%
\special{pa 1252 890}%
\special{pa 1260 920}%
\special{pa 1270 952}%
\special{pa 1278 984}%
\special{pa 1286 1016}%
\special{pa 1294 1048}%
\special{pa 1300 1082}%
\special{pa 1308 1116}%
\special{pa 1316 1150}%
\special{pa 1322 1186}%
\special{pa 1326 1220}%
\special{pa 1332 1256}%
\special{pa 1336 1290}%
\special{pa 1338 1324}%
\special{pa 1338 1358}%
\special{pa 1338 1392}%
\special{pa 1334 1426}%
\special{pa 1330 1458}%
\special{pa 1324 1488}%
\special{pa 1314 1518}%
\special{pa 1302 1546}%
\special{pa 1288 1572}%
\special{pa 1272 1598}%
\special{pa 1252 1620}%
\special{pa 1228 1642}%
\special{pa 1204 1662}%
\special{pa 1178 1682}%
\special{pa 1150 1700}%
\special{pa 1150 1700}%
\special{sp}%
%
\special{pn 20}%
\special{pa 1190 660}%
\special{pa 1160 674}%
\special{pa 1130 688}%
\special{pa 1104 704}%
\special{pa 1078 722}%
\special{pa 1058 742}%
\special{pa 1040 764}%
\special{pa 1026 790}%
\special{pa 1014 816}%
\special{pa 1004 846}%
\special{pa 998 876}%
\special{pa 992 908}%
\special{pa 988 942}%
\special{pa 986 976}%
\special{pa 986 1010}%
\special{pa 984 1046}%
\special{pa 984 1082}%
\special{pa 984 1120}%
\special{pa 982 1156}%
\special{pa 980 1192}%
\special{pa 978 1230}%
\special{pa 974 1266}%
\special{pa 972 1302}%
\special{pa 968 1336}%
\special{pa 966 1370}%
\special{pa 964 1404}%
\special{pa 964 1438}%
\special{pa 966 1468}%
\special{pa 970 1498}%
\special{pa 976 1528}%
\special{pa 986 1556}%
\special{pa 998 1582}%
\special{pa 1014 1606}%
\special{pa 1034 1628}%
\special{pa 1058 1650}%
\special{pa 1082 1670}%
\special{pa 1110 1690}%
\special{pa 1138 1708}%
\special{pa 1140 1710}%
\special{sp}%
%
\special{pn 20}%
\special{pa 3600 780}%
\special{pa 3592 812}%
\special{pa 3586 844}%
\special{pa 3578 876}%
\special{pa 3572 908}%
\special{pa 3564 938}%
\special{pa 3558 970}%
\special{pa 3554 1002}%
\special{pa 3548 1034}%
\special{pa 3546 1066}%
\special{pa 3542 1098}%
\special{pa 3540 1128}%
\special{pa 3540 1160}%
\special{pa 3540 1192}%
\special{pa 3542 1224}%
\special{pa 3546 1256}%
\special{pa 3550 1288}%
\special{pa 3554 1320}%
\special{pa 3558 1352}%
\special{pa 3564 1382}%
\special{pa 3570 1414}%
\special{pa 3576 1446}%
\special{pa 3582 1478}%
\special{pa 3588 1508}%
\special{pa 3594 1540}%
\special{pa 3600 1572}%
\special{pa 3606 1604}%
\special{pa 3612 1634}%
\special{pa 3618 1666}%
\special{pa 3624 1698}%
\special{pa 3630 1728}%
\special{pa 3636 1760}%
\special{pa 3640 1780}%
\special{sp 0.070}%
%
\special{pn 20}%
\special{pa 1860 2220}%
\special{pa 1892 2214}%
\special{pa 1924 2208}%
\special{pa 1956 2204}%
\special{pa 1988 2198}%
\special{pa 2018 2192}%
\special{pa 2050 2186}%
\special{pa 2082 2182}%
\special{pa 2114 2176}%
\special{pa 2146 2172}%
\special{pa 2178 2168}%
\special{pa 2208 2164}%
\special{pa 2240 2160}%
\special{pa 2272 2158}%
\special{pa 2304 2154}%
\special{pa 2336 2152}%
\special{pa 2368 2152}%
\special{pa 2400 2150}%
\special{pa 2432 2150}%
\special{pa 2464 2150}%
\special{pa 2496 2150}%
\special{pa 2528 2150}%
\special{pa 2560 2152}%
\special{pa 2592 2154}%
\special{pa 2624 2156}%
\special{pa 2656 2158}%
\special{pa 2686 2160}%
\special{pa 2718 2162}%
\special{pa 2750 2166}%
\special{pa 2782 2168}%
\special{pa 2814 2172}%
\special{pa 2846 2174}%
\special{pa 2878 2178}%
\special{pa 2900 2180}%
\special{sp 0.070}%
%
\special{pn 20}%
\special{pa 1800 1180}%
\special{pa 1832 1174}%
\special{pa 1864 1166}%
\special{pa 1896 1160}%
\special{pa 1928 1152}%
\special{pa 1958 1146}%
\special{pa 1990 1140}%
\special{pa 2022 1134}%
\special{pa 2054 1128}%
\special{pa 2086 1122}%
\special{pa 2116 1118}%
\special{pa 2148 1114}%
\special{pa 2180 1110}%
\special{pa 2212 1106}%
\special{pa 2244 1104}%
\special{pa 2274 1102}%
\special{pa 2306 1100}%
\special{pa 2338 1100}%
\special{pa 2370 1100}%
\special{pa 2402 1102}%
\special{pa 2432 1104}%
\special{pa 2464 1106}%
\special{pa 2496 1110}%
\special{pa 2528 1114}%
\special{pa 2560 1118}%
\special{pa 2590 1124}%
\special{pa 2622 1130}%
\special{pa 2654 1136}%
\special{pa 2686 1142}%
\special{pa 2716 1148}%
\special{pa 2748 1156}%
\special{pa 2780 1162}%
\special{pa 2812 1170}%
\special{pa 2844 1178}%
\special{pa 2874 1186}%
\special{pa 2906 1194}%
\special{pa 2938 1202}%
\special{pa 2970 1210}%
\special{pa 2970 1210}%
\special{sp 0.070}%
%
\special{pn 20}%
\special{sh 0.600}%
\special{ar 2350 690 40 40  0.0000000 6.2831853}%
%
\special{pn 20}%
\special{sh 0.600}%
\special{ar 1800 1180 40 40  0.0000000 6.2831853}%
%
\special{pn 20}%
\special{sh 0.600}%
\special{ar 2950 1200 40 40  0.0000000 6.2831853}%
%
\special{pn 20}%
\special{sh 0.600}%
\special{ar 2370 1540 40 40  0.0000000 6.2831853}%
\put(22.6000,-5.7000){\makebox(0,0)[lb]{$A_1$}}%
\put(14.8000,-12.8000){\makebox(0,0)[lb]{$A_1$}}%
\put(30.2000,-13.0000){\makebox(0,0)[lb]{$A_1$}}%
\put(23.1000,-17.8000){\makebox(0,0)[lb]{$A_1$}}%
%
\special{pn 20}%
\special{pa 2500 800}%
\special{pa 2440 870}%
\special{fp}%
%
\special{pn 20}%
\special{pa 2440 800}%
\special{pa 2502 870}%
\special{fp}%
%
\special{pn 20}%
\special{pa 2380 790}%
\special{pa 2320 860}%
\special{fp}%
%
\special{pn 20}%
\special{pa 2320 790}%
\special{pa 2382 860}%
\special{fp}%
%
\special{pn 20}%
\special{pa 2240 880}%
\special{pa 2180 950}%
\special{fp}%
%
\special{pn 20}%
\special{pa 2180 880}%
\special{pa 2242 950}%
\special{fp}%
%
\special{pn 20}%
\special{pa 2550 930}%
\special{pa 2490 1000}%
\special{fp}%
%
\special{pn 20}%
\special{pa 2490 930}%
\special{pa 2552 1000}%
\special{fp}%
%
\special{pn 20}%
\special{pa 2390 960}%
\special{pa 2330 1030}%
\special{fp}%
%
\special{pn 20}%
\special{pa 2330 960}%
\special{pa 2392 1030}%
\special{fp}%
%
\special{pn 20}%
\special{pa 2720 940}%
\special{pa 2660 1010}%
\special{fp}%
%
\special{pn 20}%
\special{pa 2660 940}%
\special{pa 2722 1010}%
\special{fp}%
%
\special{pn 20}%
\special{pa 2550 1040}%
\special{pa 2490 1110}%
\special{fp}%
%
\special{pn 20}%
\special{pa 2490 1040}%
\special{pa 2552 1110}%
\special{fp}%
%
\special{pn 20}%
\special{pa 2780 1080}%
\special{pa 2720 1150}%
\special{fp}%
%
\special{pn 20}%
\special{pa 2720 1080}%
\special{pa 2782 1150}%
\special{fp}%
%
\special{pn 20}%
\special{pa 2250 1030}%
\special{pa 2190 1100}%
\special{fp}%
%
\special{pn 20}%
\special{pa 2190 1030}%
\special{pa 2252 1100}%
\special{fp}%
%
\special{pn 20}%
\special{pa 2400 1150}%
\special{pa 2340 1220}%
\special{fp}%
%
\special{pn 20}%
\special{pa 2340 1150}%
\special{pa 2402 1220}%
\special{fp}%
%
\special{pn 20}%
\special{pa 2080 970}%
\special{pa 2020 1040}%
\special{fp}%
%
\special{pn 20}%
\special{pa 2020 970}%
\special{pa 2082 1040}%
\special{fp}%
%
\special{pn 20}%
\special{pa 2660 1180}%
\special{pa 2600 1250}%
\special{fp}%
%
\special{pn 20}%
\special{pa 2600 1180}%
\special{pa 2662 1250}%
\special{fp}%
%
\special{pn 20}%
\special{pa 2810 1220}%
\special{pa 2750 1290}%
\special{fp}%
%
\special{pn 20}%
\special{pa 2750 1220}%
\special{pa 2812 1290}%
\special{fp}%
%
\special{pn 20}%
\special{pa 1970 1070}%
\special{pa 1910 1140}%
\special{fp}%
%
\special{pn 20}%
\special{pa 1910 1070}%
\special{pa 1972 1140}%
\special{fp}%
%
\special{pn 20}%
\special{pa 2090 1140}%
\special{pa 2030 1210}%
\special{fp}%
%
\special{pn 20}%
\special{pa 2030 1140}%
\special{pa 2092 1210}%
\special{fp}%
%
\special{pn 20}%
\special{pa 2260 1210}%
\special{pa 2200 1280}%
\special{fp}%
%
\special{pn 20}%
\special{pa 2200 1210}%
\special{pa 2262 1280}%
\special{fp}%
%
\special{pn 20}%
\special{pa 2110 1270}%
\special{pa 2050 1340}%
\special{fp}%
%
\special{pn 20}%
\special{pa 2050 1270}%
\special{pa 2112 1340}%
\special{fp}%
%
\special{pn 20}%
\special{pa 2250 1360}%
\special{pa 2190 1430}%
\special{fp}%
%
\special{pn 20}%
\special{pa 2190 1360}%
\special{pa 2252 1430}%
\special{fp}%
%
\special{pn 20}%
\special{pa 2520 1270}%
\special{pa 2460 1340}%
\special{fp}%
%
\special{pn 20}%
\special{pa 2460 1270}%
\special{pa 2522 1340}%
\special{fp}%
%
\special{pn 20}%
\special{pa 2660 1310}%
\special{pa 2600 1380}%
\special{fp}%
%
\special{pn 20}%
\special{pa 2600 1310}%
\special{pa 2662 1380}%
\special{fp}%
%
\special{pn 20}%
\special{pa 2400 1300}%
\special{pa 2340 1370}%
\special{fp}%
%
\special{pn 20}%
\special{pa 2340 1300}%
\special{pa 2402 1370}%
\special{fp}%
%
\special{pn 20}%
\special{pa 2490 1400}%
\special{pa 2430 1470}%
\special{fp}%
%
\special{pn 20}%
\special{pa 2430 1400}%
\special{pa 2492 1470}%
\special{fp}%
%
\special{pn 4}%
\special{pa 1330 1390}%
\special{pa 990 1050}%
\special{fp}%
\special{pa 1330 1330}%
\special{pa 990 990}%
\special{fp}%
\special{pa 1320 1260}%
\special{pa 1000 940}%
\special{fp}%
\special{pa 1310 1190}%
\special{pa 1010 890}%
\special{fp}%
\special{pa 1300 1120}%
\special{pa 1020 840}%
\special{fp}%
\special{pa 1270 1030}%
\special{pa 1040 800}%
\special{fp}%
\special{pa 1260 960}%
\special{pa 1060 760}%
\special{fp}%
\special{pa 1230 870}%
\special{pa 1090 730}%
\special{fp}%
\special{pa 1210 790}%
\special{pa 1130 710}%
\special{fp}%
\special{pa 1320 1440}%
\special{pa 990 1110}%
\special{fp}%
\special{pa 1310 1490}%
\special{pa 990 1170}%
\special{fp}%
\special{pa 1290 1530}%
\special{pa 990 1230}%
\special{fp}%
\special{pa 1270 1570}%
\special{pa 980 1280}%
\special{fp}%
\special{pa 1250 1610}%
\special{pa 980 1340}%
\special{fp}%
\special{pa 1210 1630}%
\special{pa 970 1390}%
\special{fp}%
\special{pa 1180 1660}%
\special{pa 970 1450}%
\special{fp}%
\special{pa 1150 1690}%
\special{pa 990 1530}%
\special{fp}%
%
\special{pn 4}%
\special{pa 2480 160}%
\special{pa 2330 310}%
\special{fp}%
\special{pa 2430 150}%
\special{pa 2270 310}%
\special{fp}%
\special{pa 2370 150}%
\special{pa 2230 290}%
\special{fp}%
\special{pa 2310 150}%
\special{pa 2170 290}%
\special{fp}%
\special{pa 2240 160}%
\special{pa 2130 270}%
\special{fp}%
\special{pa 2170 170}%
\special{pa 2080 260}%
\special{fp}%
\special{pa 2070 210}%
\special{pa 2030 250}%
\special{fp}%
\special{pa 2530 170}%
\special{pa 2390 310}%
\special{fp}%
\special{pa 2580 180}%
\special{pa 2460 300}%
\special{fp}%
\special{pa 2630 190}%
\special{pa 2540 280}%
\special{fp}%
\special{pa 2680 200}%
\special{pa 2630 250}%
\special{fp}%
\put(9.9000,-21.0000){\makebox(0,0)[lb]{dense}}%
%
\special{pn 8}%
\special{pa 1400 1940}%
\special{pa 2040 1450}%
\special{fp}%
\special{sh 1}%
\special{pa 2040 1450}%
\special{pa 1976 1476}%
\special{pa 1998 1482}%
\special{pa 2000 1506}%
\special{pa 2040 1450}%
\special{fp}%
\put(31.9000,-21.5000){\makebox(0,0)[lb]{four ends tend to infinity}}%
%
\special{pn 8}%
\special{pa 3150 2160}%
\special{pa 2960 2190}%
\special{fp}%
\special{sh 1}%
\special{pa 2960 2190}%
\special{pa 3030 2200}%
\special{pa 3014 2182}%
\special{pa 3024 2160}%
\special{pa 2960 2190}%
\special{fp}%
%
\special{pn 8}%
\special{pa 4020 1970}%
\special{pa 3720 1730}%
\special{fp}%
\special{sh 1}%
\special{pa 3720 1730}%
\special{pa 3760 1788}%
\special{pa 3762 1764}%
\special{pa 3786 1756}%
\special{pa 3720 1730}%
\special{fp}%
\put(29.3000,-6.2000){\makebox(0,0)[lb]{$\Sol(\th)_{\bR}$}}%
\end{picture}%
\end{center}
\caption{Real Cayley cubic $\Sol(\th)_{\bR}$ with four 
$A_1$-singularities}
\label{fig:cayley}
\end{figure}
\begin{example}[Stratum of type $\mbox{\bm $A_1^{\oplus3}$}$] 
\label{ex:A13}
This is the locus discussed by Dubrovin and Mazzocco \cite{DM}, although 
they made use of a different parametrization of the character variety. 
On this stratum they showed that there are exactly five algebraic 
solutions up to some equivalence. 
\end{example}
\begin{example}[Stratum of type $\mbox{\boldmath $A_1^{\oplus2}$}$] 
\label{ex:A12}
This stratum is not well understood yet. 
We content ourselves with giving an example, 
the orbit in Figure \ref{fig:A12}. 
It is a finite $G$-orbit of degree $6$ with parameters 
$\th = (2\sqrt{2},2\sqrt{2},3,4) \in \Th$, 
which is the $\rh$-image of $\k = (1/4,0,0,1/12,5/12) \in \K$, 
certainly a point of type $A_1^{\oplus2}$. 
This $G$-orbit is also a $G(2)$-orbit of degree $6$. 
\end{example}
\section{Tetrahedral Theorem} \label{sec:tetra} 
The strata belonging to the sequence (S2) admit a unified treatment. 
\begin{theorem} \label{thm:tetra}
On any $F_4^{(1)}$-stratum belonging to the sequence $(\mathrm{S}2)$, 
there is no non-Riccati algebraic solutions of degree $d \ge 7$ without 
univalent local branches at any fixed singular point. 
\end{theorem}
This theorem can be used to classify all algebraic solutions on the 
strata belonging to the sequence (S2). 
We may refer to Theorem \ref{thm:tetra} as the 
{\sl Tetrahedral Theorem} for the following reasons. 
\begin{table}[t] 
\begin{center}
\begin{tabular}{|c|c|}
\hline
\vspace{-3mm} & \\
$D_4^{(1)}$-strata along sequence (S2) & skeletons of tetrahedron \\[1mm]
\hline
\vspace{-3mm} & \\
one stratum of abstract type $A_1$ & one $3$-cell \\[1mm]
$\downarrow$ & $\downarrow$ \\[1mm]
four strata of abstract type $A_2$ & four faces \\[1mm]
$\downarrow$ & $\downarrow$ \\[1mm]
six strata of abstract type $A_3$ & six edges \\[1mm]
$\downarrow$ & $\downarrow$ \\[1mm]
four strata of abstract type $D_4$ & four vertices \\[1mm]
\hline
\end{tabular}
\end{center}
\caption{A parallelism in adjacency relations}
\label{tab:parallel}
\end{table}
\begin{remark}[Parallelism] \label{rem:tetra}
There is a parallelism as in Table \ref{tab:parallel} between the 
adjacency relations for the $D_4^{(1)}$-strata along the sequence 
(S2) and that for the skeletons of the (regular) tetrahedron. 
This parallelism is not by chance. 
Behind it there exists an interesting story starting with the 
algebraic geometry of Painlev\'e VI and ending up with some 
elementary geometry of a regular tetrahedron of edge length $\sqrt{2}$. 
Indeed, in the course of establishing Theorem \ref{thm:tetra} we 
come across the regular tetrahedron in 
Figure \ref{fig:tetrahedron} (right), which lies in the $3$-dimensional 
space with coordinates $(m_0/d, m_1/d, m_{\infty}/d)$, where $d$ is the 
degree of the algebraic solution under consideration and 
$(m_0,m_1,m_{\infty})$ is a triplet of positive integers encoding 
certain information of how the algebraic solution branches 
at the fixed singular points 
$z = 0, 1, \infty$. 
A detailed explanation can be found in \cite{Iwasaki2}. 
\end{remark}
\par
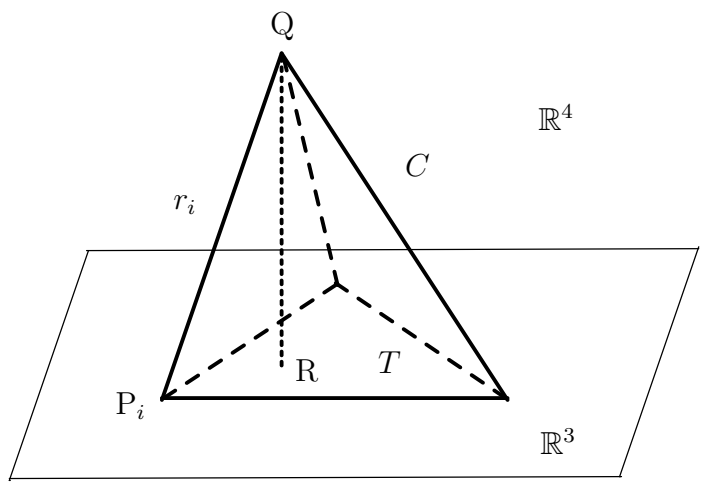
\begin{figure}[t]
\begin{center}
\unitlength 0.1in
\begin{picture}( 35.7000, 24.7000)(  6.2000,-26.5000)
%
\special{pn 8}%
\special{pa 620 2640}%
\special{pa 3790 2630}%
\special{fp}%
%
\special{pn 8}%
\special{pa 1020 1450}%
\special{pa 4190 1440}%
\special{fp}%
%
\special{pn 8}%
\special{pa 1030 1450}%
\special{pa 620 2650}%
\special{fp}%
\special{pa 620 2650}%
\special{pa 620 2650}%
\special{fp}%
%
\special{pn 8}%
\special{pa 4190 1430}%
\special{pa 3780 2630}%
\special{fp}%
\special{pa 3780 2630}%
\special{pa 3780 2630}%
\special{fp}%
%
\special{pn 20}%
\special{pa 1410 2220}%
\special{pa 3200 2220}%
\special{fp}%
%
\special{pn 20}%
\special{pa 2030 420}%
\special{pa 1410 2230}%
\special{fp}%
%
\special{pn 20}%
\special{pa 2030 420}%
\special{pa 3200 2230}%
\special{fp}%
%
\special{pn 20}%
\special{pa 2310 1630}%
\special{pa 1420 2220}%
\special{da 0.070}%
%
\special{pn 20}%
\special{pa 2310 1620}%
\special{pa 3190 2220}%
\special{da 0.070}%
%
\special{pn 20}%
\special{pa 2040 440}%
\special{pa 2320 1630}%
\special{da 0.070}%
%
\special{pn 20}%
\special{pa 2030 440}%
\special{pa 2030 2050}%
\special{dt 0.054}%
\put(19.7000,-3.5000){\makebox(0,0)[lb]{$\rQ$}}%
\put(11.7000,-23.3000){\makebox(0,0)[lb]{$\rP_i$}}%
\put(14.7000,-12.6000){\makebox(0,0)[lb]{$r_i$}}%
\put(21.0000,-21.3000){\makebox(0,0)[lb]{$\rR$}}%
\put(25.3000,-20.9000){\makebox(0,0)[lb]{$T$}}%
\put(33.7000,-25.1000){\makebox(0,0)[lb]{$\bR^3$}}%
\put(33.6000,-8.2000){\makebox(0,0)[lb]{$\bR^4$}}%
\put(26.7000,-10.6000){\makebox(0,0)[lb]{$C$}}%
\end{picture}%
\end{center}
\caption{$4$-dimensional cone $C$ over the tetrahedron $T$}
\label{fig:geometry}
\end{figure}

We explain what kind of elementary geometry comes up.  
Let $T = \rP_1\rP_2\rP_3\rP_4 \subset \bR^3$ be a regular 
tetrahedron with edge length $\sqrt{2}$ as in 
Figure \ref{fig:tetrahedron} (left); 
$C = \rQ\rP_1\rP_2\rP_3\rP_4 \subset \bR^4$ the cone over 
the base $T$ with side lengths $\ol{\rQ\rP}_i = r_i$ for 
$i = 1,2,3,4$, as in Figure \ref{fig:geometry}; and let 
$\rR$ be the orthogonal projection of the vertex $\rQ$ down 
to the $3$-space $\bR^3$ that contains the tetrahedron $T$. 
Moreover let $\ora{\rR}$ and $\ora{\rP_i}$ denote the 
position vectors of the points $\rR$ and $\rP_i$ respectively. 
Write 
\[
\ora{\rR} = 
\a_1 \ora{\rP_1} + \a_2 \ora{\rP_2} + \a_3 \ora{\rP_3} 
+ \a_4 \ora{\rP_4}, 
\]
in terms of the barycentric coordinates 
$\a = (\a_1,\a_2,\a_3,\a_4) \in \bR^4$ where 
$\a_1 + \a_2 +\a_3 + \a_4 = 1$. 
In the Painlev\'e situation, $T$ is the tetrahedron of 
Figure \ref{fig:tetrahedron} (right) and the vertices 
$\rP_i$ $(i = 1,2,3,4)$ are just those of the latter tetrahedron. 
A basic lemma we need is the following. 
\begin{lemma} \label{lem:geometry} 
If the side lengths $r_i$ $(i = 1,2,3,4)$ are chosen as 
\begin{equation} \label{eqn:radii}
\left\{
\begin{array}{rcl}
r_1^2 &=& (\k_1-1)^2+\k_2^2+\k_3^2+\k_4^2, \\[1mm] 
r_2^2 &=& \k_1^2+(\k_2-1)^2+\k_3^2+\k_4^2, \\[1mm]
r_3^2 &=& \k_1^2+\k_2^2+(\k_3-1)^2+\k_4^2, \\[1mm] 
r_4^2 &=& \k_1^2+\k_2^2+\k_3^2+(\k_4-1)^2,  
\end{array}
\right.
\end{equation}
with $\k = (\k_0,\k_1,\k_2,\k_3,\k_4) \in \K_{\bR}$, then 
\begin{equation} \label{eqn:territory}
\ol{\rQ\rR}^2 = \k_0^2, \qquad \a_i = \k_i + \frac{\k_0}{2} 
\qquad (i = 1,2,3,4). 
\end{equation}
\end{lemma}
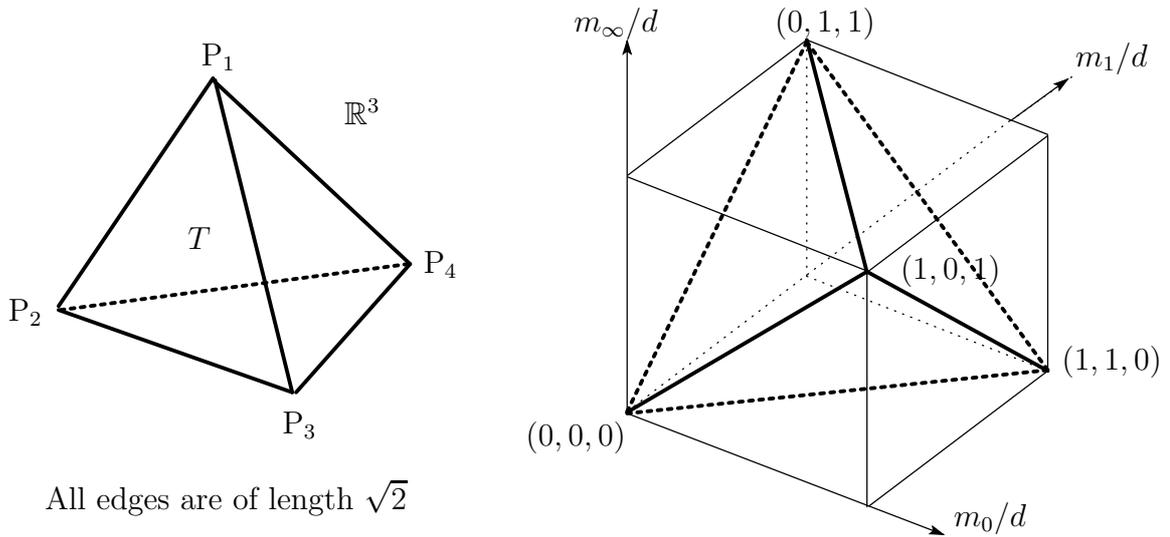
\begin{figure}[t]
\begin{center}
\unitlength 0.1in
\begin{picture}( 55.2600, 27.5200)(  3.5000,-31.9200)
%
\special{pn 8}%
\special{pa 4806 3044}%
\special{pa 5732 2338}%
\special{fp}%
%
\special{pn 8}%
\special{pa 3558 1316}%
\special{pa 4484 610}%
\special{fp}%
%
\special{pn 8}%
\special{pa 3558 1326}%
\special{pa 4806 1820}%
\special{fp}%
%
\special{pn 8}%
\special{pa 4806 1816}%
\special{pa 5732 1110}%
\special{fp}%
%
\special{pn 8}%
\special{pa 4484 610}%
\special{pa 5732 1104}%
\special{fp}%
%
\special{pn 8}%
\special{pa 4796 1820}%
\special{pa 4796 3048}%
\special{fp}%
%
\special{pn 8}%
\special{pa 5732 1114}%
\special{pa 5732 2344}%
\special{fp}%
%
\special{pn 8}%
\special{pa 4484 616}%
\special{pa 4484 1844}%
\special{dt 0.045}%
%
\special{pn 8}%
\special{pa 4484 1844}%
\special{pa 5732 2338}%
\special{dt 0.045}%
%
\special{pn 8}%
\special{pa 3554 2550}%
\special{pa 5536 1032}%
\special{dt 0.045}%
%
\special{pn 8}%
\special{pa 5540 1032}%
\special{pa 5858 802}%
\special{dt 0.045}%
%
\special{pn 8}%
\special{pa 5540 1032}%
\special{pa 5834 816}%
\special{fp}%
\special{sh 1}%
\special{pa 5834 816}%
\special{pa 5768 840}%
\special{pa 5790 848}%
\special{pa 5792 872}%
\special{pa 5834 816}%
\special{fp}%
%
\special{pn 20}%
\special{pa 4484 610}%
\special{pa 4796 1824}%
\special{fp}%
%
\special{pn 20}%
\special{pa 4796 1820}%
\special{pa 3554 2554}%
\special{fp}%
%
\special{pn 20}%
\special{pa 4796 1820}%
\special{pa 5732 2338}%
\special{fp}%
%
\special{pn 20}%
\special{pa 3554 2560}%
\special{pa 5732 2334}%
\special{dt 0.054}%
%
\special{pn 20}%
\special{pa 4484 620}%
\special{pa 3554 2550}%
\special{dt 0.054}%
%
\special{pn 20}%
\special{pa 4484 620}%
\special{pa 5732 2334}%
\special{dt 0.054}%
%
\special{pn 8}%
\special{pa 3554 2554}%
\special{pa 3554 616}%
\special{fp}%
\special{sh 1}%
\special{pa 3554 616}%
\special{pa 3534 682}%
\special{pa 3554 668}%
\special{pa 3574 682}%
\special{pa 3554 616}%
\special{fp}%
%
\special{pn 8}%
\special{pa 3554 2560}%
\special{pa 5190 3192}%
\special{fp}%
\special{sh 1}%
\special{pa 5190 3192}%
\special{pa 5136 3150}%
\special{pa 5140 3174}%
\special{pa 5122 3188}%
\special{pa 5190 3192}%
\special{fp}%
\put(30.3000,-27.7000){\makebox(0,0)[lb]{$(0,0,0)$}}%
\put(58.0400,-23.7200){\makebox(0,0)[lb]{$(1,1,0)$}}%
\put(43.2000,-6.1000){\makebox(0,0)[lb]{$(0,1,1)$}}%
\put(49.7000,-19.0000){\makebox(0,0)[lb]{$(1,0,1)$}}%
\put(52.4700,-31.9200){\makebox(0,0)[lb]{$m_{0}/d$}}%
\put(58.7600,-7.8800){\makebox(0,0)[lb]{$m_{1}/d$}}%
\put(32.8000,-6.1000){\makebox(0,0)[lb]{$m_{\infty}/d$}}%
%
\special{pn 20}%
\special{pa 1410 810}%
\special{pa 610 2000}%
\special{fp}%
%
\special{pn 20}%
\special{pa 1420 820}%
\special{pa 1820 2450}%
\special{fp}%
%
\special{pn 20}%
\special{pa 1410 820}%
\special{pa 2430 1780}%
\special{fp}%
%
\special{pn 20}%
\special{pa 610 2020}%
\special{pa 1820 2450}%
\special{fp}%
%
\special{pn 20}%
\special{pa 2430 1780}%
\special{pa 1840 2450}%
\special{fp}%
%
\special{pn 20}%
\special{pa 620 2020}%
\special{pa 2420 1780}%
\special{dt 0.054}%
\put(17.7000,-26.9000){\makebox(0,0)[lb]{$\rP_3$}}%
\put(25.0000,-18.5000){\makebox(0,0)[lb]{$\rP_4$}}%
\put(13.5000,-7.7000){\makebox(0,0)[lb]{$\rP_1$}}%
\put(3.5000,-21.0000){\makebox(0,0)[lb]{$\rP_2$}}%
\put(12.8000,-17.0000){\makebox(0,0)[lb]{$T$}}%
\put(5.4000,-30.9000){\makebox(0,0)[lb]{All edges are of length $\sqrt{2}$}}%
\put(20.9000,-10.5000){\makebox(0,0)[lb]{$\bR^3$}}%
\end{picture}%
\end{center}
\vspace{-2mm}
\caption{Tetrahedron for the Tetrahedral Theorem}
\label{fig:tetrahedron}
\end{figure}
\par 
It is difficult to explain in short words why the choice 
(\ref{eqn:radii}) is natural in our situation and we refer 
to \cite{Iwasaki2} for a detailed explanation. 
Anyway, in the course of establishing Theorem \ref{thm:tetra} 
we encounter a sort of {\sl territory problem}, where the 
territory of the vertex $\rP_i$ is the $3$-dimensional open 
ball $B_i := B(\rP_i,r_i)$ of radius $r_i$ with center at 
the point $\rP_i$. 
To explain what this problem is all about, we begin by  
stating a key observation as in the following lemma. 
\begin{lemma} \label{lem:territory}
If $\PVI(\k)$ with $\k \in \K_{\bR}$ admits a non-Riccati algebraic 
solution without univalent local branches at any fixed singular point, 
then the balls $B_i$ $(i = 1,2,3,4)$ must have at least one points 
in common. 
\end{lemma}
As the contraposition of this lemma, if the four balls 
have no points in common then there is no algebraic solution with 
the prescribed property. 
Now a natural question is when they have points in common 
and when not. 
Let us restrict our attention to the case where the point $\rR$ 
lies in the interior of $T$, that is, where the barycentric 
coordinates $\a = (\a_1,\a_2,\a_3,\a_4)$ satisfy the inequalities
\begin{equation} \label{eqn:bary}
\a_i > 0 \qquad (i = 1,2,3,4). 
\end{equation}
In this case, if the balls $B_i$ $(i = 1,2,3,4)$ have at least one 
points in common, then $\rR$ must be such a point in common. 
With this observation we are in a position to give the following. 
\par\medskip\noindent
{\bf Sketch of the proof of Theorem \ref{thm:tetra}}. 
The proof is by contradiction. 
Assume that $\PVI(\k)$ has a non-Riccati algebraic solution of 
degree $d \ge 7$ without univalent local branches at any fixed singular 
point. 
Then we must have $\k \in \K_{\bR}$ from Theorem \ref{thm:rationality}. 
After applying a suitable B\"acklund transformation we may 
assume that $\k$ lies in the (closed) fundamental $W(D_4^{(1)})$-alcove 
$\{ \k \in \K_{\bR}\,:\, \k_i \ge 0 \,\, (i = 0,1,2,3,4) \}$. 
Now assume that $\k$ lies on the stratum of type $A_1$. 
Then there is a unique index $i_0 \in \{0,1,2,3,4\}$ such that 
$\k_{i_0} = 0$ and $\k_i > 0$ for the remaining indices $i$. 
After applying a further B\"acklund transformation we may assume that 
$i_0 = 0$, namely, that $\k_0 = 0$ and $\k_i > 0$ for $i = 1,2,3,4$.  
So it follows from formula (\ref{eqn:territory}) that 
\begin{equation} \label{eqn:territory2}
\ol{\rQ\rR} = 0, \qquad \a_i = \k_i > 0  \qquad (i = 1,2,3,4). 
\end{equation}
The former condition in (\ref{eqn:territory2}) means that $\rR = \rQ$ 
and hence $\ol{\rR\rP}_i = \ol{\rQ\rP}_i = r_i$ so that $\rR$ is a 
point of the boundary sphere $\partial B_i$. 
Since $B_i$ is an open ball, $\rR$ does not belong to $B_i$. 
On the other hand the latter condition in (\ref{eqn:territory2}) 
means that condition (\ref{eqn:bary}) is satisfied so that 
$\rR$ must belong to $B_i$, a contradiction. 
Similar arguments are feasible on the other strata of 
the sequence (S2). \hfill $\Box$ 
\section{The Big Open} \label{sec:BigOpen}
On the big open $\K_{\emptyset} = \K - \Wall(D_4^{(1)})$ we are still
 distant from the complete classification, but we are already able 
to confine all finite orbits into a rather thin subset of the 
real character variety $\Sol(\th)_{\bR}$ (see \cite{Iwasaki2}). 
In dealing with this stratum it is necessary to distinguish the 
two subsets $\Wall(D_4^{(1)})$ and $\Wall(F_4^{(1)})$ of the parameter 
space $\K$, where the former is the union of all reflecting hyperplanes 
for the reflection group $W(D_4^{(1)})$ and the latter is its 
counterpart for the group $W(F_4^{(1)})$. 
Note that there is the strict inclusion 
$\Wall(D_4^{(1)}) \subset \Wall(F_4^{(1)})$. 
In the parameter level almost all algebraic solutions on this stratum 
seem to exist on the set $\Wall(F_4^{(1)})-\Wall(D_4^{(1)})$. 
In fact, Boalch's ``generic" icosahedral solution \cite{Boalch2} 
(see also item (2) of Remark \ref{rem:rationality}) is the 
only instance outside $\Wall(F_4^{(1)})$ known so far 
(as of September 8, 2008). 

\end{document}